# Notes on Theory of Quadratic Residues

*by Roberto Volpe*[*]

***Abstract:***

The Law of Quadratic Reciprocity was conjectured by Euler and Legendre who both found an incomplete proof. Gauss called this law "*Theorema Fundamentale*", and he was the first who gave a complete proof, he also highlighted the equivalence of his formulation with those of Euler and Legrendre.
Hereby notes gives a overview of the Theory of Quadratic Residues using a classical approach with some application to Diophantine Equations, such as Two Square Theorem and Pythagorean Quadruplets.

[*] Mathematics Enthusiastic from Vercelli (Italy)

# Note sulla Teoria dei Residui Quadratici

*di Roberto Volpe (*)*

(*) Vercellese Cultore della Matematica

# *Indice*

# 1. Introduzione.

*"Se $p$ è un numero primo della forma $4n+1$ , $+p$ sarà un residuo o non residuo di un qualsiasi numero primo che preso positivamente è un residuo o un non residuo di $p$. Se $p$ è della forma $4n+3$, $-p$ avrà la stessa proprietà."*

Carl Friedrich Gauss DISQUISITIONES ARITMETICAE - Articolo 131

La legge della reciprocità quadratica fu congetturata e parzialmente dimostrata da Eulero e da Legendre, Gauss chiamava tale legge nella formulazione sopra riportata *"Theorema Fundamentale"*, e fu il primo a darne una dimostrazione soddisfacente evidenziando l'equivalenza con le formulazioni di Eulero e Legendre.

Queste note non hanno la pretesa di dire nulla di nuovo sull'argomento, ma sono semplicemente degli appunti di studio dell'autore. In teoria un lettore interessato al termine della lettura di queste note dovrebbe avere una panoramica generale di tutta la teoria sottesa dalla soluzione generale della congruenza:

$$aX^2 + bX + c \equiv 0 \pmod{n}$$

Mi è sembrato naturale aggiungere anche il famoso teorema di Fermat sugli interi somma di due quadrati e i suoi ampliamenti, in quanto derivabili dalla teoria dei residui quadratici, e conseguentemente parlare degli interi di Gauss dandone una applicazione alla risoluzione della equazione diofantea:

$$X^2 + Y^2 = cZ^2$$

della equazione Diofantea:

$$X^2 + Y^2 = Z^l$$

e infine della equazione Diofantea:

$$X^2 + Y^2 + Z^2 = W^2$$

ottenendo così una formula per le quadruple pitagoriche.
Analizzando queste equazioni diofantee si generalizza in tre direzioni differenti quanto noto sulle terne pitagoriche, e infatti in tutti e tre i casi le ritroveremo come caso particolare.

Sarebbe naturale proseguire parlando della teoria delle forme quadratiche binarie e ternarie a coefficienti interi e quindi parlare della equazione Diofantea:

$$aX^2 + bY^2 = cZ^2$$

nota in letteratura come equazione di Legendre. Per ora non sono stati aggiunti questi argomenti alla presente nota, preferendo limitare la trattazione ad un breve cenno per quanto strettamente necessario a spiegare la dimostrazione di Legendre delle legge della reciprocità quadratica.

Si assumono come noti alcuni concetti base di Teoria di numeri e di Algebra, e per essi si può fare riferimento ai testi [A] e [MA] indicati in Bibliografia o alla vastissima bibliografia presente on-line.

## 2. *Residui quadratici modulo un numero primo*

Iniziamo con l'occuparci della risolubilità della congruenza $X^2 \equiv a \pmod{p}$ con $p$ primo dispari.

*Definizione 2.1.*

Sia $p$ un primo dispari ed $a$ un intero non divisibile per $p$. Diciamo che $a$ è un residuo quadratico di $p$ quando la $X^2 \equiv a \pmod{p}$ ha almeno una soluzione, altrimenti $a$ è un residuo non quadratico.

*Osservazione 2.2*

Per ogni intero $a$ non divisibile per $2$ (cioè ogni intero dispari) la congruenza $X^2 \equiv a \pmod{2}$ ha soluzioni, infatti il quadrato di un qualsiasi numero dispari è dispari e quindi congruo modulo $2$ ad $a$. Il numero primo $2$ costituisce una eccezione rispetto agli altri numeri primi perché possiede una sola classe di soluzioni tra loro tutte congrue modulo $2$. Per gli altri primi invece esistono sempre due classi di soluzioni non congrue come dimostrato nella proposizione seguente. Questo giustifica il fatto che $2$ non sia contemplato nella definizione 2.1. In ogni caso vedremo meglio nel paragrafo 7 che per le potenze di 2, come $4$ e $2^e$ con $e \geq 3$, è necessaria una trattazione dettagliata.

*Proposizione 2.3*

Se $a$ è un residuo quadratico del numero primo dispari $p$, allora la congruenza $X^2 \equiv a \pmod{p}$ ha esattamente due soluzioni non congruenti modulo $p$ date da $X = \pm b$ con $b^2 \equiv a \pmod{p}$.

Dimostrazione

Se $b^2 \equiv a \pmod{p}$ allora $b$ e $-b$ sono soluzioni e sono incongruenti modulo $p$. Infatti se per assurdo $p \mid 2b$ allora siccome $MCD(2, p) = 1$, segue che $p \mid b$ da cui essendo $b^2 \equiv a \pmod{p}$ segue che $p \mid a$ in contraddizione con la definizione 3.1.1. Se $c$ è una soluzione di $X^2 \equiv a \pmod{p}$ allora $p \mid c^2 - a$ ed essendo $c^2 - a = c^2 - b^2 + b^2 - a$ segue che $p \mid c^2 - b^2$, cioè risulta che $p \mid c - b$ oppure $p \mid c + b$, ovvero $c \equiv b \pmod{p}$ oppure $c \equiv -b \pmod{p}$. □

*Proposizione 2.4*

I residui quadratici modulo $p$ si ottengono calcolando $b^2$ modulo $p$ per il solo intervallo $1 \leq b \leq (p-1)/2$. Gli altri quadrati per $(p-1)/2 < b < p-1$ si ripetono in ordine inverso, pertanto esattamente metà dei residui non nulli (modulo $p$) sono residui quadratici (modulo $p$).

Dimostrazione

Basta osservare che $(p-b)^2 \equiv b^2 \pmod{p}$, quindi i residui che si ottengono facendo variare $1 \leq b \leq (p-1)/2$ si ripetono in senso inverso rispetto a quelli ottenuti facendo variare $b$ nell'intervallo $(p-1)/2 < b < p-1$. Non ci possono essere due residui quadratici congruenti generati da $1 \leq b' < b'' \leq (p-1)/2$ o da $(p-1)/2 < (p-b') < (p-b'') \leq p-1$ perché essendo $b'' - b' < p$ se $p \mid b'' - b'$ segue che $b'' = b'$, quindi i residui quadratici sono esattamente metà dei residui quadratici non nulli. □

## 3. *Il simbolo di Legendre e il criterio di Eulero.*

*Definizione 3.1*

Il simbolo di Legendre $\left(\frac{a}{p}\right)$ è definito come segue:

$\left(\frac{a}{p}\right) = +1$ se $p \nmid a$ e $x^2 \equiv a \pmod p$ ha soluzioni.

$\left(\frac{a}{p}\right) = -1$ se $p \nmid a$ e $x^2 \equiv a \pmod p$ non ha soluzioni.

$\left(\frac{a}{p}\right) = 0$ se $p \mid a$.

*Proposizione 3.2*

Sia $p$ un primo dispari:

(a) $a \equiv b \pmod p \Rightarrow \left(\frac{a}{p}\right) = \left(\frac{b}{p}\right)$.

(b) ***Criterio di Eulero***: $\left(\frac{a}{p}\right) \equiv a^{(p-1)/2} \pmod p$.

(c) $\left(\frac{-1}{p}\right) = (-1)^{(p-1)/2} = \begin{cases} +1 & \text{se } p \equiv 1 \pmod 4 \\ -1 & \text{se } p \equiv 3 \pmod 4 \end{cases}$.

(d) $\left(\frac{ab}{p}\right) = \left(\frac{a}{p}\right)\left(\frac{b}{p}\right)$.

Dimostrazione

(a) Se $a \equiv b \pmod p$ allora se $x^2 \equiv a \pmod p$ allora anche $x^2 \equiv b \pmod p$ e viceversa, perché $x^2 - a = x^2 - b + b - a$, quindi $\left(\frac{a}{p}\right) = \left(\frac{b}{p}\right)$.

(b) Sia $\mathbb{F}_p = \mathbb{Z}/p\mathbb{Z}$ il campo delle classi di resto modulo $p$ e $\mathbb{F}_p^* = \{[1]_p, \ldots, [p-1]_p\}$ il relativo gruppo moltiplicativo associato (vedi [A] oppure [MA]).
Consideriamo $x_1 \in \{1, \ldots, p-1\}$ allora $[x_1] \cdot \mathbb{F}_p^* = \{[1]_p, \ldots, [p-1]_p\} = \{[x_1], \ldots, [x_1(p-1)]\} = \mathbb{F}_p^*$, per cui uno ed uno solo dei numeri $x_1$, $2x_1, \ldots, (p-1)x_1$ è congruente ad $a$. Chiamiamo $x_1'$ l'unico numero per cui $x_1 x_1' \equiv a \pmod p$. Fissato $a$ con $p \nmid a$ consideriamo separatamente i due casi che possono presentarsi:

- $\left(\frac{a}{p}\right) = 1$

  In tal caso esistono esattamente due soluzioni di $X^2 \equiv a \pmod p$, cioè, $c^2 \equiv a \pmod p$ e $(p-c)^2 \equiv a \pmod p$, quindi $c(p-c) \equiv -a \pmod p$ i rimanenti $p-3$ elementi di $\{1,2,\ldots,p-1\} - \{c, p-c\}$ sono costituiti da coppie di interi distinti $(x_1, x_1')$ tali che $x_1 x_1' \equiv a \pmod p$.
  Segue che $1 \cdot 2 \cdot \ldots \cdot (p-1) = (p-1)!$ commutando i fattori del prodotto risulta congruente al numero $-a \cdot a^{\frac{p-3}{2}} = -a^{\frac{p-1}{2}}$.
  In sintesi: $(p-1)! \equiv -a^{(p-1)/2} \pmod p$

- $\left(\frac{a}{p}\right) = -1$

  In tal caso esitono $(p-1)/2$ coppie $(x_1, x_1')$ tali che $x_1 x_1' \equiv a \pmod p$. Quindi in tal caso si ha la congruenza:
  $(p-1)! \equiv a^{(p-1)/2} \pmod p$.

Come caso particolare abbiamo il ***teorema di Wilson*** per $a = 1$, infatti $x^2 \equiv 1 \pmod p$ ha come soluzioni $x = \pm 1$ quindi abbiamo:
$(p-1)! \equiv -1 \pmod p$
Infine ricordando la definizione di simbolo di Legendre abbiamo la formula cercata:
$\left(\frac{a}{p}\right) \equiv a^{(p-1)/2} \pmod p$.

(c) Come applicazione del punto (b) abbiamo per $a = -1$ che:
$\left(\frac{-1}{p}\right) \equiv (-1)^{(p-1)/2} \pmod p$ quindi
Se $(p-1)/2$ è pari cioè $p \equiv 1 \pmod 4$ allora $-1$ è un residuo quadratico modulo $p$.
Se $(p-1)/2$ è dispari cioè $p \equiv 3 \pmod 4$ allora $-1$ non è un residuo quadratico modulo $p$.

(d) Ecco una ulteriore applicazione di (b):
$\left(\frac{ab}{p}\right) \equiv (ab)^{(p-1)/2} \equiv a^{(p-1)/2} b^{(p-1)/2} \equiv \left(\frac{a}{p}\right)\left(\frac{b}{p}\right) \pmod p$
siccome $-2 \le \left(\frac{ab}{p}\right) - \left(\frac{a}{p}\right)\left(\frac{b}{p}\right) \le 2$ e $p > 2$ segue che $\left(\frac{ab}{p}\right) = \left(\frac{a}{p}\right)\left(\frac{b}{p}\right)$. □

*Corollario 3.3*

Se $p \equiv 1 \pmod 4$ allora $\left(\frac{-a}{p}\right) = \left(\frac{a}{p}\right)$ per tutti gli $a$.

Se $p \equiv 3 \pmod 4$ allora $\left(\frac{-a}{p}\right) = -\left(\frac{a}{p}\right)$ per tutti gli $a$.

Dimostrazione
Se $p \equiv 1 \pmod 4$ allora $(p-1)/2$ è pari quindi applicando il criterio di Eulero dato che $(-a)^{(p-1)/2} = a^{(p-1)/2}$ segue il risultato.
Se invece $p \equiv 3 \pmod 4$ allora $(p-1)/2$ è dispari quindi applicando il criterio di Eulero dato che $(-a)^{(p-1)/2} = -a^{(p-1)/2}$ segue il risultato. □

La risposta alla domanda "quali primi $p$ hanno $2$ come residuo quadratico" ci porta ad uno dei più importanti risultati della teoria elementare dei numeri. Sperimentalmente per i primi piccoli possiamo congetturare la risposta "i primi congruenti $\pm 1 \pmod 8$"

$\left(\frac{2}{p}\right) = 1$ per $p = 7,17,23,31,41,47,71...$

$\left(\frac{2}{p}\right) = -1$ per $p = 3,5,11,13,19,29,37,43...$

Più in generale $\left(\frac{a}{p}\right)$, per $a$ fissato, dipende dal residuo di $p$ modulo $4a$, e questa è proprio una delle forme della legge della reciprocità quadratica di Gauss, inizialmente congetturata da Eulero.

## 4. *La Legge di Reciprocità Quadratica (LRQ).*

*Proposizione 4.1-**Lemma di Gauss***

Sia $p$ un primo dispari e $a$ un intero non divisibile per $p$. Si consideri il minimo residuo dell'intero $ka$ per $1 \le k \le (p-1)/2$ ridotto in modo da essere compreso nell'intervallo $[-p/2, p/2]$. Se il numero di questi residui che sono negativi è $s$, allora:

$$\left(\frac{a}{p}\right) = (-1)^s.$$

Dimostrazione

Abbiamo che per applicazione del *criterio di Eulero*:

$$(1 \cdot a) \cdot (2 \cdot a) \ldots (\frac{(p-1)}{2} \cdot a) = \frac{(p-1)}{2}!\, a^{(p-1)/2} = \frac{(p-1)}{2}! \left(\frac{a}{p}\right)$$

inoltre risulta:

$$(1 \cdot a) \cdot (2 \cdot a) \ldots (\frac{(p-1)}{2} \cdot a) = (-1)^s \frac{(p-1)}{2}!$$ infatti:

- Se $ka \equiv k'a \pmod p$, dato che $p \nmid a$, allora $p \mid (k'-k)$ da cui essendo $-\frac{p-1}{2} \le k'-k \le \frac{p-1}{2}$ segue che $k = k'$.
- Se per assurdo $ka \equiv -k'a \pmod p$, dato che $p \nmid a$, allora $p \mid (k'+k)$ da cui essendo $0 < k'+k < p$ segue l'assurdo $k = -k'$ dato che $kk' > 0$.
- Non può risultare $ka \equiv 0 \pmod p$ dato che $p \nmid a$ e $p \nmid k$.

Quindi ciascuno dei numeri $ka$ per $1 \le k \le (p-1)/2$ è congruo ad uno ed uno solo dei numeri dei numeri non nulli del sistema di residui modulo $p$: $A = \left\{-\frac{(p-1)}{2}, \ldots, -2, -1, 0, 1, 2, \ldots \frac{(p-1)}{2}\right\}$, e numeri differenti sono congruenti a residui differenti in valore assoluto. Per cui $(1 \cdot a) \cdot (2 \cdot a) \ldots (\frac{(p-1)}{2} \cdot a)$ ha come valore assoluto $\frac{(p-1)}{2}!$, e come segno il prodotto dei segni dei residui negativi, cioè $(-1)^s$.
Per semplice uguaglianza delle due formule sopra semplificando il fattore comune si ottiene il risultato. □

Siamo ora in grado di rispondere alla nostra domanda relativa a primi dispari che hanno come residuo quadratico $2$.

*Proposizione 4.2*

Sia $p$ un primo dispari allora:

$$\left(\frac{2}{p}\right) = (-1)^{(p^2-1)/8} = \begin{cases} +1 & se \quad p \equiv \pm 1 \pmod 8 \\ -1 & se \quad p \equiv \pm 3 \pmod 8 \end{cases}.$$

Dimostrazione

Si tratta di una applicazione del Lemma di Gauss, contiamo il numero $s$ di residui in $\left[-\frac{p}{2}, \frac{p}{2}\right]$ dei numeri $2$, $4$, ..., $p-1$ con segno negativo. Condizione necessaria e sufficiente a generare un residuo negativo e che $2k > \frac{p-1}{2}$, cioè $k > \frac{p-1}{4}$.

- Se $\frac{p-1}{4}$ è intero si hanno esattamente $\frac{p-1}{2} - \frac{p-1}{4} = \frac{p-1}{4}$ residue negativi, quindi:

- Se $\frac{p-1}{4}$ è dispari, cioè $p \equiv -3 \pmod 8$ e $\left(\frac{2}{p}\right) = -1$.
- Se $\frac{p-1}{4}$ è pari, cioè $p \equiv 1 \pmod 8$ e $\left(\frac{2}{p}\right) = 1$.

➢ Se $\frac{p-1}{4}$ non è intero, allora l'ultimo residuo positivo è dato da $2k = \frac{p-3}{2}$ e quindi si hanno esattamente $\frac{p-1}{2} - \frac{p-3}{4} = \frac{p+1}{4}$ residui negativi, quindi:

- Se $\frac{p+1}{4}$ è dispari, cioè $p \equiv 3 \pmod 8$ e $\left(\frac{2}{p}\right) = -1$.
- Se $\frac{p+1}{4}$ è pari, cioè $p \equiv -1 \pmod 8$ e $\left(\frac{2}{p}\right) = 1$.

Notando che $\frac{(p^2-1)}{8} = \frac{(p+1)}{2}\frac{(p-1)}{4} = \frac{(p-1)}{2}\frac{(p+1)}{4}$ ha la stessa parità di $s$ si ha formula sintetica. □

*Teorema 4.3* -***Reciprocità Quadratica***

Siano $p$ e $q$ primi dispari distinti, allora:

$$\left(\frac{p}{q}\right)\left(\frac{q}{p}\right) = (-1)^{\frac{(p-1)}{2}\frac{(q-1)}{2}}$$

Cioè:

Se $p \equiv 1 \pmod 4$ o $q \equiv 1 \pmod 4$, allora $\left(\frac{q}{p}\right) = \left(\frac{p}{q}\right)$

Se $p \equiv 3 \pmod 4$ e $q \equiv 3 \pmod 4$, allora $\left(\frac{q}{p}\right) = -\left(\frac{p}{q}\right)$

Dimostrazione

Consideriamo i seguenti due insiemi:

$$S = \left\{1,\dots,\frac{p-1}{2}\right\}, T = \left\{1,\dots,\frac{q-1}{2}\right\}$$

Per il lemma di Gauss:

$$\left(\frac{q}{p}\right) = (-1)^{s_1} \text{ dove: } s_1 = \left|\left\{h \in S \mid \exists h' \in S : qh \in [-h']_p\right\}\right|$$

$$\left(\frac{p}{q}\right) = (-1)^{s_2} \text{ dove: } s_2 = \left|\left\{k \in T \mid \exists k' \in T : pk \in [-k']_q\right\}\right|$$

Quindi $\left(\frac{p}{q}\right)\left(\frac{q}{p}\right) = (-1)^{s_1+s_2}$, per cui si tratta di trovare la parità di $s_1 + s_2$.

Se $qh \equiv -h' \pmod p$ con $h \in S$, allora esiste unico $k \in \mathbb{Z}$ tale che $pk - qh = -h' \in S$, cioè, $0 < pk - qh \le \frac{(p-1)}{2}$. Il numero $k$ deve stare in $T$, infatti:

$$pk > qh > -h' > 0 \text{ e } pk \le \frac{p-1}{2} + qh \le (q+1)\frac{p-1}{2} < p\frac{q+1}{2}$$

quindi $0 < k \le \frac{q-1}{2} < \frac{q+1}{2}$, cioè $k \in T$. Segue che:

$$s_1 = |U_1|, \text{ con } U_1 = \left\{(h,k) \in S \times T \mid 0 < pk - qh \le \frac{(p-1)}{2}\right\}$$

In maniera perfettamente analoga si dimostra che:

$s_2 = |U_2|$, con $U_2 = \left\{(h,k) \in S \times T \mid -\frac{(q-1)}{2} \leq pk - qh < 0\right\}$

Siccome non ci sono coppie $(h,k) \in S \times T$ tali che $pk = qh$, ne segue che:

$U = \left\{(h,k) \in S \times T \mid -\frac{(q-1)}{2} \leq pk - qh \leq \frac{(p-1)}{2}\right\} = U_1 \cup U_2$, inoltre siccome $U_1 \cap U_2 = \emptyset$

$s_1 + s_2 = |U| = \left|\left\{(h,k) \in S \times T \mid -\frac{(q-1)}{2} \leq pk - qh \leq \frac{(p-1)}{2}\right\}\right|$

$U$ è invariante rispetto alla rotazione $(h,k) \xrightarrow{\vartheta} \left(\frac{p+1}{2} - h, \frac{q+1}{2} - k\right)$ di $180^0$, rispetto al centro $\left(\frac{p+1}{4}, \frac{q+1}{4}\right)$ del rettangolo $R = \left[0, \frac{p+1}{2}\right] \times \left[0, \frac{q+1}{2}\right]$, cioè $\vartheta(U) = U$, infatti:

$pk' - qh' = p\left(\frac{q+1}{2} - k\right) - q\left(\frac{p+1}{2} - h\right) = \frac{p-q}{2} - (pk - qh)$

quindi, $-\frac{(q-1)}{2} \leq pk - qh \leq \frac{(p-1)}{2} \Leftrightarrow -\frac{(q-1)}{2} \leq pk' - qh' \leq \frac{(p-1)}{2}$.

La rotazione $\vartheta$ ha il solo punto fisso $\left(\frac{p+1}{4}, \frac{q+1}{4}\right)$, quindi

$\left(\frac{q}{p}\right)\left(\frac{p}{q}\right) = -1 \Leftrightarrow |U|$ è dispari $\Leftrightarrow \frac{p+1}{4} \in \mathbb{Z} \wedge \frac{q+1}{4} \in \mathbb{Z} \Leftrightarrow \left(p \equiv 3 \pmod 4\right) \wedge \left(q \equiv 3 \pmod 4\right)$. □

## 5. *Congettura di Eulero e dimostrazione alternativa della LRQ.*

La legge della reciprocità quadratica prima di essere dimostrata da Gauss fu congetturata da Eulero nella seguente forma.

*Teorema 5.1*-***Congettura di Eulero***
Siano $p$ e $q$ due primi dispari e $a$ un intero coprimo con $pq$.

**(i)** Se $p \equiv q \pmod{4a}$, allora $\left(\frac{a}{p}\right) = \left(\frac{a}{q}\right)$.

**(ii)** Se $p \equiv -q \pmod{4a}$, allora $\left(\frac{a}{p}\right) = \left(\frac{a}{q}\right)$.

Come Gauss scrive nell'articolo 151 di [DA], Eulero dimostrò in maniera lacunosa e quindi in realtà congetturò solo, la seguente proposizione:

*"esistono numeri $r'$, $r''$, $r'''$, etc., $n'$, $n''$, $n'''$, etc., $< 4a$ tali che i divisori primi di $x^2 - a$ sono contenuti in una delle forme $r' + k(4a)$, $r'' + k(4a)$, $r''' + k(4a)$, etc., e tutti i primi non divisori nella forma $n' + k(4a)$, $n'' + k(4a)$, $n''' + k(4a)$, etc., con $k$ variabile sugli interi."*

Sempre nell'articolo 151 di [DA] Gauss fa notare la proposizione sopra può essere facilmente dimostrata dalla seguente proposizione dimostrata da Eulero nello studio dell'equazione: $fx^2 + gy^2 = hz^2$:

*"se l'equazione è risolubile per un valore $h = s$, sarà anche risolubile per ogni altro valore congruente a $s$ relativamente al modulo $4fg$ a condizione che sia un numero primo"*

Infatti se $a = -fg$ e $fX^2 + gY^2 = sZ^2$ con $s$ primo è risolubile, allora $T^2 \equiv a \pmod{s}$ ha soluzioni, cioè $\left(\frac{a}{s}\right) = 1$, infatti se $fx_0{}^2 + gy_0{}^2 = sz_0{}^2$, allora: $\left(\frac{a}{s}\right) = \left(\frac{ay_0{}^2}{s}\right) = \left(\frac{f^2x_0{}^2 - sfz_0{}^2}{s}\right) = 1$.

Quindi dato che ogni primo della forma $h = s + k(4a)$ rende risolubile $fx^2 + gy^2 = hz^2$, segue che anche che per ognuno di tali primi si ha $\left(\frac{a}{h}\right) = 1$. Quello che in generale è più difficile da dimostrare che se $\left(\frac{a}{s}\right) = 1$ allora esistono $f$ e $g$ tali che $a = -fg$ e $fX^2 + gY^2 = sZ^2$ è risolubile, darò maggiori dettagli nell'osservazione 7.9 dimostrando il caso in cui $a = q$ primo.

Dopo questa breve digressione storica che dimostra come Eulero sia andato vicino alla dimostrazione del teorema 5.1, vediamo come il teorema 3.4.1 sia completamente equivalente alla legge della reciprocità quadratica. Prima di tutto dimostriamo che la legge della reciprocità quadratica implica il teorema 5.1.

*Lemma 5.2*
Sia $r$ un primo dispari e $p$, $q$ due primi dispari distinti da $r$, allora valgono le seguenti:

**(i)** Se $p \equiv q \pmod{4r}$, allora $\left(\frac{r}{p}\right) = \left(\frac{r}{q}\right)$.

**(ii)** Se $p \equiv -q \pmod{4r}$, allora $\left(\frac{r}{p}\right) = \left(\frac{r}{q}\right)$.

Dimostrazione
Applicando la legge della reciprocità quadratica

$\left(\frac{r}{p}\right)\left(\frac{p}{r}\right) = (-1)^{\frac{r-1}{2}\frac{p-1}{2}}$, $\left(\frac{r}{q}\right)\left(\frac{q}{r}\right) = (-1)^{\frac{r-1}{2}\frac{q-1}{2}}$;

Nel caso **(i)**. posto $p = q + 4kr$ abbiamo dividendo membro a membro che:

$$\frac{\left(\frac{r}{p}\right)}{\left(\frac{r}{q}\right)} \cdot \frac{\left(\frac{p}{r}\right)}{\left(\frac{q}{r}\right)} = (-1)^{\frac{r-1}{2}\frac{p-q}{2}} = \left((-1)^{\frac{r-1}{2}kr}\right)^2 = 1,$$

dato che $\left(\frac{p}{r}\right) = \left(\frac{q+4kr}{r}\right) = \left(\frac{q}{r}\right)$, si ha il risultato cercato.

Nel caso **(ii)** posto $p = -q + 4kr$ abbiamo moltiplicando membro a membro:

$$\left(\frac{r}{p}\right)\left(\frac{r}{q}\right)\left(\frac{p}{r}\right)\left(\frac{q}{r}\right) = (-1)^{\frac{r-1}{2}\left(\frac{p+q}{2}+1\right)} = (-1)^{\frac{r-1}{2}},$$

dato che $\left(\frac{p}{r}\right) = \left(\frac{-q+4kr}{r}\right) = \left(\frac{-q}{r}\right) = \left(\frac{-1}{r}\right)\left(\frac{q}{r}\right) = (-1)^{\frac{r-1}{2}}\left(\frac{q}{r}\right)$, si ha di nuovo il risultato cercato. □

*Dimostrazione del Teorema 5.1*

Sia $a = r_1^{e_1}...r_n^{e_n}$ la scomposizione in fattori primi di $a$ segue che:

$$\left(\frac{a}{p}\right) = \left(\frac{r_1^{e_1}...r_n^{e_n}}{p}\right) = \left(\frac{r_1}{p}\right)^{e_1}...\left(\frac{r_n}{p}\right)^{e_n}, \left(\frac{a}{q}\right) = \left(\frac{r_1^{e_1}...r_n^{e_n}}{q}\right) = \left(\frac{r_1}{q}\right)^{e_1}...\left(\frac{r_n}{q}\right)^{e_n},$$

siccome se $p \equiv \pm q \pmod{4a}$, allora $p \equiv \pm q \pmod{4r_j}$ per ogni $1 \le j \le n$, quindi per il lemma 3.4.2, abbiamo che: $\left(\frac{r_j}{p}\right) = \left(\frac{r_j}{q}\right)$ per ogni $1 \le j \le n$, da cui segue il risultato cercato. □

*Teorema 5.3*

Siano $p$ e $q$ due primi dispari distinti:

P1. $a$ un intero coprimo con $pq$.

**(iii)** Se $p \equiv q \pmod{4a}$, allora $\left(\frac{a}{p}\right) = \left(\frac{a}{q}\right)$.

**(iv)** Se $p \equiv -q \pmod{4a}$, allora $\left(\frac{a}{p}\right) = \left(\frac{a}{q}\right)$.

Se e solo se

P2. $\left(\frac{q}{p}\right) = \left(\frac{p}{q}\right)$ se $p \equiv 1 \pmod 4$ o $q \equiv 1 \pmod 4$,

$\left(\frac{q}{p}\right) = -\left(\frac{p}{q}\right)$ se $p \equiv 3 \pmod 4$ e $q \equiv 3 \pmod 4$.

Dimostrazione

Abbiamo gia visto con il lemma 5.2 e con la dimostrazione del teorema 5.1. che P2 implica P1. Verifichiamo ora che P1 implica P2.

Se $p \equiv q \pmod 4$, allora $p - q = 4a$ e $p \nmid a \vee q \nmid a$ dato che $p$ e $q$ sono primi distinti. Quindi

$$\left(\frac{q}{p}\right) = \left(\frac{p-4a}{p}\right) = \left(\frac{-4a}{p}\right) = \left(\frac{-1}{p}\right)\left(\frac{a}{p}\right) = (-1)^{\frac{p-1}{2}}\left(\frac{a}{q}\right) = (-1)^{\frac{p-1}{2}}\left(\frac{4a}{q}\right) = (-1)^{\frac{p-1}{2}}\left(\frac{p-q}{q}\right) = (-1)^{\frac{p-1}{2}}\left(\frac{p}{q}\right)$$

da cui:

$\left(\frac{q}{p}\right) = \left(\frac{p}{q}\right)$ se $p \equiv q \equiv 1 \pmod 4$, e $\left(\frac{q}{p}\right) = -\left(\frac{p}{q}\right)$ se $p \equiv q \equiv 3 \pmod 4$

Se $p \equiv -q \pmod 4$, allora $p + q = 4a$ e $p \nmid a \vee q \nmid a$ dato che $p$ e $q$ sono primi distinti. Quindi:

$$\left(\frac{q}{p}\right)=\left(\frac{-p+4a}{p}\right)=\left(\frac{4a}{p}\right)=\left(\frac{a}{p}\right)=\left(\frac{a}{q}\right)=\left(\frac{4a}{q}\right)=\left(\frac{p+q}{q}\right)=\left(\frac{p}{q}\right)$$

da cui:

$\left(\frac{q}{p}\right)=\left(\frac{p}{q}\right)$ se $p\equiv -q\equiv 1 \pmod 4$, oppure $-p\equiv q\equiv 1 \pmod 4$. □

Il teorema 5.1 può essere dimostrato direttamente usando il Lemma di Gauss, e in questo modo abbiamo una dimostrazione alternativa della legge della reciprocità quadratica.

*Dimostrazione del teorema 3.4.1 usando il Lemma di Gauss*

Per il lemma di Gauss sappiamo che $\left(\frac{a}{p}\right)=(-1)^{s_1}$ dove, posto $S=\left\{1,\dots,\frac{p-1}{2}\right\}$, si definisce $s_1=\left|\left\{h\in S \mid \exists h'\in S : ah\in[-h']_p\right\}\right|$, cioè il numero dei residui dei numeri $aS=\left\{a,2a,\dots,\frac{p-1}{2}a\right\}$ che cadono nell'intervallo $\left]-\frac{p}{2},0\right[$. Equivalentemente siccome ogni numero che è congruo modulo $p$ ad un numero nell'intervallo $\left]-\frac{p}{2},0\right[$, è anche congruo modulo $p$ ad un numero nell'intervallo $\left]\frac{p}{2},p\right[$, $s$ è anche il numero di residui positivi modulo $p$ dei numeri di $aS$ che cade nell'intervallo $\left]\frac{p}{2},p\right[$.

Siccome due numeri di $aS$ non possono essere congruenti modulo $p$, segue che $s_1$ è anche il numero di elementi di $aS$ che cade in uno dei seguenti intervalli:

$I_1=\left]\frac{p}{2},p\right[$, $I_2=\left]\frac{3}{2}p,2p\right[$ ,…, $I_k=\left]\frac{2k-1}{2}p,kp\right[$ con $k=\left\lfloor\frac{a}{2}\right\rfloor$.

> Apriamo una breve parentesi e verifichiamo perché $k=\left\lfloor\frac{a}{2}\right\rfloor$.
>
> Se $a=2m$, allora $\frac{p-1}{2}2m=(p-1)m<pm=p\left\lfloor\frac{a}{2}\right\rfloor=pk$.
>
> Se $a=2m+1$, allora $\frac{p-1}{2}(2m+1)=(p-1)\left(m+\frac{1}{2}\right)<p\left(m+\frac{1}{2}\right)=p\left\lfloor\frac{a}{2}\right\rfloor+\frac{p}{2}=p\left(k+\frac{1}{2}\right)$,
>
> siccome nell'intervallo $\left]kp,\left(k+\frac{1}{2}\right)p\right[$ non possono cadere numeri con minimi residui positivi in $\left]\frac{p}{2},p\right[$, l'ultimo in cui può caderne uno è $I_k=\left]\frac{2k-1}{2}p,kp\right[$.

Abbiamo quindi posto $I=\bigcup_{j=1}^{k} I_j$ che risulta $s_1=|aS\cap I|$.

Gli interi di $aS$ sono tutti multipli di $a$, e ci può essere al più un intero per ogni intervallo $I_j$, e viceversa se $ha\in I$, allora $ha\in aS$.

Quindi il numero di interi contenuti in $\frac{1}{a}I=\bigcup_{j=1}^{k}\frac{1}{a}I_j=\bigcup_{j=1}^{k}\left]\frac{2j-1}{2a}p,\frac{j}{a}p\right[$ coincide con $s$, cioè

$s_1=|aS\cap I|=\left|\mathbb{Z}\cap\frac{1}{a}I\right|$.

Un identico ragionamento può essere ripetuto per $q$, dove poniamo: $T = \left\{1,\dots,\frac{q-1}{2}\right\}$,

$J_i = \left]\frac{2k-1}{2}q, kq\right[$ con $i = 1,\dots,k$; $J = \bigcup_{i=1}^{k} J_i$, e per il lemma di Gauss $\left(\frac{a}{q}\right) = (-1)^{s_2}$ dove ora risulta $s_2 = |aT \cap I| = \left|\mathbb{Z} \cap \frac{1}{a}J\right|$.

Posto ora che sia $p = q + 4ah$ e consideriamo la generica coppia di intervalli $I_j, J_j$ risulta che:

$$|I_j \cap \mathbb{Z}| = \left|\left]\frac{2j-1}{2}p, \frac{j}{a}p\right[ \cap \mathbb{Z}\right| = \left|\left]\frac{2j-1}{2a}q + 2h, \frac{j}{a}q + 4h\right[ \cap \mathbb{Z}\right| =$$

$$\left|\left]\frac{2j-1}{2a}q + 2h, \frac{j}{a}q + 2h\right[ \cap \mathbb{Z}\right| + \left|\left]\frac{j}{a}q + 2h, \frac{j}{a}q + 4h\right[ \cap \mathbb{Z}\right| = |J_j \cap \mathbb{Z}| + 2h$$

quindi $|I_j \cap \mathbb{Z}| \equiv |I_j \cap \mathbb{Z}| \pmod 2$ da cui infine $s_1 \equiv s_2 \pmod 2$.

Nel caso invece che sia $p = -q + 4ah$ consideriamo di nuovo la generica coppia di intervalli:

$$|I_j \cap \mathbb{Z}| = \left|\left]\frac{2j-1}{2a}p, \frac{j}{a}p\right[ \cap \mathbb{Z}\right| = \left|\left]-\frac{j}{a}p, -\frac{2j-1}{2a}p\right[ \cap \mathbb{Z}\right| = \left|\left]\frac{j}{a}q - 4h, \frac{2j-1}{2a}q - 2h,\right[ \cap \mathbb{Z}\right| =$$

$$\left|\left]\frac{j}{a}q - 4h, \frac{j}{a}q - 2h\right[ \cap \mathbb{Z}\right| - \left|\left]\frac{2j-1}{2a}q - 2h, \frac{j}{a}q - 2h\right[ \cap \mathbb{Z}\right| = 2h - |J_j \cap \mathbb{Z}|$$

quindi $|I_j \cap \mathbb{Z}| \equiv -|I_j \cap \mathbb{Z}| \pmod 2$ da cui infine $s_1 \equiv -s_2 \pmod 2$, da cui anche $s_1 \equiv s_2 \pmod 2$ dato che se la somma di due numeri è pari lo anche la loro differenza e viceversa. □

La seguente proposizione consente di definire in dettaglio dato un primo $q$, quali primi dispari hanno $q$ come residuo quadratico.

*Proposizione 5.4*
Sia $q$ un primo dispari; $q$ è un residuo quadratico modulo il primo dispari $p \neq q$ se e solo se $p$ è congruente modulo $4q$ a uno dei seguenti interi: $\pm 1^2$, $\pm 3^2$, $\pm 5^2$, ..., $\pm (q-2)^2$.

Dimostrazione
Se $p \equiv (2a+1)^2 \pmod{4q}$ allora in particolare $p \equiv 1 \pmod 4$ applicando la legge della reciprocità quadratica:

$$\left(\frac{q}{p}\right) = \left(\frac{p}{q}\right)(-1)^{\frac{p-1}{2}\frac{q-1}{2}} = \left(\frac{(2a+1)^2}{q}\right)(-1)^{\frac{p-1}{2}\frac{q-1}{2}} = (-1)^{\frac{p-1}{2}\frac{q-1}{2}} = 1$$

Se $p \equiv -(2a+1)^2 \pmod{4q}$ allora $p \equiv -1 \pmod 4$ e di nuovo applicando la legge della reciprocità quadratica:

$$\left(\frac{q}{p}\right) = \left(\frac{p}{q}\right)(-1)^{\frac{p-1}{2}\frac{q-1}{2}} = \left(\frac{-(2a+1)^2}{q}\right)(-1)^{\frac{p-1}{2}\frac{q-1}{2}} = \left(\frac{-1}{q}\right)^{\frac{p-1}{2}\frac{q-1}{2}} = (-1)^{\frac{p+1}{2}\frac{q-1}{2}} = 1$$

Viceversa, se $\left(\frac{q}{p}\right) = 1$, per il criterio di Eulero:

$$\left((-1)^{\frac{p-1}{2}}\right)^{\frac{q-1}{2}} \equiv \left(\frac{(-1)^{\frac{p-1}{2}}}{q}\right) \pmod q$$

da cui segue avendo a che fare con numeri in valore assoluto uguali ad 1 che:

$$\left((-1)^{\frac{p-1}{2}}\right)^{\frac{q-1}{2}} = \left(\frac{(-1)^{\frac{p-1}{2}}}{q}\right)$$

Applicando la legge della reciprocità quadratica

$\left(\frac{p}{q}\right) = (-1)^{\frac{p-1}{2}\frac{q-1}{2}} = \left(\frac{(-1)^{\frac{p-1}{2}}}{q}\right)$ quindi $\left(\frac{p(-1)^{\frac{p-1}{2}}}{q}\right) = \left(\frac{p}{q}\right)^2 = 1$, cioè esiste $0 < x < q$ tale che:

$$x^2 \equiv p(-1)^{\frac{p-1}{2}} \pmod{q}$$

siccome $x^2 \equiv (q-x)^2 \pmod{q}$ e uno dei due numeri tra $x$ e $q-x$ è dispari si può assumere $x$ dispari, quindi $x^2 \equiv 1 \pmod{4}$.

Se $p \equiv 1 \pmod{4}$ allora $x^2 \equiv p \pmod{q}$ da cui essendo $p$ dispari segue che $x^2 \equiv p \pmod{4q}$.

Se $p \equiv -1 \pmod{4}$ allora $x^2 \equiv -p \pmod{q}$ da cui $x^2 \equiv -p \pmod{4q}$. □

Quindi ad esempio la regola che definisce i primi per cui $3$ è un residuo quadratico è la seguente:

$$\left(\frac{3}{p}\right) = 1 \Leftrightarrow \begin{cases} p \equiv 1 \pmod{12} \\ p \equiv -1 \pmod{12} \end{cases}$$

## 6. *Il calcolo del simbolo di Legendre e il simbolo di Jacobi*

Il simbolo di Legendre è completamente moltiplicativo, cioè $\left(\frac{ab}{p}\right)=\left(\frac{a}{p}\right)\left(\frac{b}{p}\right)$ e conosciamo come calcolare i valori di $\left(\frac{-1}{p}\right)$ e $\left(\frac{2}{p}\right)$, infine la legge di reciprocità quadratica ci dice che $\left(\frac{q}{p}\right)$ è strettamente legato a $\left(\frac{p}{q}\right)$, un fatto tutt'altro che ovvio. Quindi possediamo un metodo molto efficace per il calcolo di $\left(\frac{a}{p}\right)$, con $p$ primo dispari. Infatti, se $a$ è un intero non nullo e non divisibile per $p$:

*Passo 1*

Se $a < p$ passare al Passo 2. Se $a > p$, applicando l'algoritmo della divisione è sempre possibile ricondurci al caso $a < p$, infatti si scrive $a = a_1 p + a_0$ con $0 < a_0 < p$ e risulta

$$\left(\frac{a}{p}\right)=\left(\frac{a_1 p + a_0}{p}\right)=\left(\frac{a_0}{p}\right).$$

*Passo 2*

Usando il teorema fondamentale dell'aritmetica scriviamo $a=(-1)^k 2^{f_0} q_1^{f_1} q_2^{f_2} \ldots q_r^{f_r}$ con $k \in \{0,1\}$, $f_i \geq 0$ e ogni primo dispari $q_i$ è distinto da $p$.

Usando il fatto che il simbolo di Legendre è completamente moltiplicativo abbiamo che:

$$\left(\frac{a}{p}\right)=\left(\frac{-1}{p}\right)^k\left(\frac{2}{p}\right)^{f_0}\left(\frac{q_1}{p}\right)^{f_2}\ldots\left(\frac{q_r}{p}\right)^{f_r} \tag{6.1}$$

*Passo 3*

Applicando la Proposizione 3.2(c) si calcola $\left(\frac{-1}{p}\right)$.

*Passo 4*

Applicando la Proposizione 4.2 si calcola $\left(\frac{2}{p}\right)$.

*Passo 5*

Per il calcolo $\left(\frac{q}{p}\right)$, con $q$ primo dispari, se $q > p$ si riparte dal Passo 1, se invece $q < p$ si applica la reciprocità quadratica e quindi si riparte dal Passo 1.

Lo svantaggio di questo metodo e che dobbiamo eseguire una fattorizzazione di un intero ad ogni stadio intermedio necessario per il calcolo del simbolo di Legendre nel caso di $q$ primo dispari. Per superare questa difficoltà generalizziamo il simbolo di Legendre in modo da considerare un intero arbitrario invece di un primo dispari.

### *Definizione 6.1-***Il Simbolo di Jacobi**

Sia $a \in \mathbb{Z}$, e sia $n$ un intero dispari arbitrario con fattorizzazione $n = \pm p_1^{e_1} p_2^{e_2} \ldots p_s^{e_s}$ allora definiamo il simbolo di Jacobi come segue:

$$\left(\frac{a}{n}\right)=\left(\frac{a}{p_1}\right)^{e_1}\left(\frac{a}{p_2}\right)^{e_2}\cdots\left(\frac{a}{p_s}\right)^{e_s} \qquad (6.2)$$

Che gode delle seguenti proprietà che estendono quelle del simbolo di Legendre:

*Proposizione 6.2*

**(i)** $\left(\frac{a}{n}\right)=0$ se e solo se $MCD(a,n)>1$.

**(ii)** $\left(\frac{a}{n}\right)=\left(\frac{a}{|n|}\right)$

**(iii)** Se $a \equiv b \pmod{n}$, allora $\left(\frac{a}{n}\right)=\left(\frac{b}{n}\right)$.

**(iv)** $\left(\frac{ab}{n}\right)=\left(\frac{a}{n}\right)\left(\frac{b}{n}\right)$

**(v)** $\left(\frac{a^2}{n}\right)=1$ se $MCD(a,n)=1$.

Dimostrazione
Derivano immediatamente dalla formula (6.2). □

Comunque la cosa più importante e che anche il simbolo di Jacobi obbedisce alla legge della reciprocità quadratica.

*Teorema 6.3*
Siano $m$ ed $n$ due interi positivi e dispari, allora abbiamo che:

**(i)** $\left(\frac{-1}{n}\right)=(-1)^{\frac{n-1}{2}}$

**(ii)** $\left(\frac{2}{n}\right)=(-1)^{\frac{n^2-1}{8}}$

**(iii)** $\left(\frac{m}{n}\right)=(-1)^{\frac{m-1}{2}\frac{n-1}{2}}\left(\frac{n}{m}\right)$

Dimostrazione
La dimostrazione in tutti e tre i casi si basa sul fatto che le applicazioni di $\mathbb{N}$ in $\{-1,1\}$ definite dalle formule $n \mapsto (-1)^{\frac{n-1}{2}}$ e $n \mapsto (-1)^{\frac{n^2-1}{8}}$ sono moltiplicative, e l'applicazione di $\mathbb{N}\times\mathbb{N}$ in $\{-1,1\}$ definita dalla formula $(m,n)\mapsto(-1)^{\frac{m-1}{2}\frac{n-1}{2}}$ è bi-moltiplicativa. Siano $n = p_1^{e_1}p_2^{e_2}\ldots p_s^{e_s}$ e $m = q_1^{f_1}q_2^{f_2}\ldots q_r^{f_r}$.

**(i)** $n \equiv 1 \pmod{4}$ se e solo se $\sum_{j=1}^{s} e_j\frac{(p_j-1)}{2} \equiv \frac{(n-1)}{2} \pmod{2}$.

Infatti essendo $n$ dispari può essere il prodotto di numeri primi congrui ad $1 \pmod{4}$ o congrui a $3 \pmod{4}$ e $n \equiv 1 \pmod{4}$ se e solo se è pari il numero di primi congrui $3 \pmod{4}$ con $e$ dispari. Quindi:

$$\left(\frac{-1}{n}\right)=\left(\frac{-1}{p_1}\right)^{e_1}\left(\frac{-1}{p_2}\right)^{e_2}\cdots\left(\frac{-1}{p_r}\right)^{e_r}=(-1)^{\sum_{j=1}^{s}e_j\frac{(p_j-1)}{2}}=(-1)^{\frac{(n-1)}{2}}$$

**(ii)** $n \equiv \pm 1 \pmod{8}$ se e solo se $\sum_{j=1}^{s} e_j\frac{(p_j^2-1)}{8} \equiv \frac{(n^2-1)}{8} \pmod{2}$.

Infatti essendo $n$ dispari può solo essere il prodotto di numeri primi congrui ad $\pm 1 \pmod 8$ o ad $\pm 3 \pmod 8$, e $n \equiv \pm 1 \pmod 8$ se e solo se è pari il numero di primi congrui a $\pm 3 \pmod 8$ con $e$ dispari. Quindi:

$$\left(\frac{2}{n}\right) = \left(\frac{2}{p_1}\right)^{e_1}\left(\frac{2}{p_2}\right)^{e_2}\cdots\left(\frac{2}{p_r}\right)^{e_r} = (-1)^{\sum_{j=1}^{r} e_j \frac{(p_j^2-1)}{8}} = (-1)^{\frac{(n^2-1)}{8}}.$$

**(iii)** In tal caso abbiamo che:

$$\left(\frac{m}{n}\right) = \prod_{j=1}^{r}\left(\frac{m}{p_j}\right) = \prod_{j=1}^{s}\prod_{i=1}^{r}\left(\frac{q_i}{p_j}\right)^{e_j f_i} = \prod_{j=1}^{s}\prod_{i=1}^{r}\left(\frac{p_j}{q_i}\right)^{e_j f_i}(-1)^{e_j f_i \frac{(p_j-1)}{2}\frac{(q_i-1)}{2}} =$$

$$= \left(\prod_{i=1}^{r}\left(\frac{n}{q_i}\right)^{f_i}\right)(-1)^{\sum_{j=1}^{s} e_j \frac{(p_j-1)}{2}\sum_{j=1}^{r} f_i \frac{(q_i-1)}{2}} = \left(\frac{n}{m}\right)(-1)^{\sum_{j=1}^{s} e_j \frac{(p_j-1)}{2}\sum_{j=1}^{r} f_i \frac{(q_i-1)}{2}}$$

Siccome risulta per quanto detto in **(i)**:

$\sum_{j=1}^{r} f_j \frac{(q_i-1)}{2} \equiv \frac{(m-1)}{2} \pmod 2$ abbiamo che $\sum_{j=1}^{r} f_j \frac{(q_j-1)}{2} = \frac{(m-1)}{2} + 2k$ con $k \in \mathbb{Z}$, e

$\sum_{i=1}^{s} e_j \frac{(p_j-1)}{2} \equiv \frac{(n-1)}{2} \pmod 2$ abbiamo che $\sum_{j=1}^{s} e_j \frac{(p_j-1)}{2} = \frac{(n-1)}{2} + 2h$ con $h \in \mathbb{Z}$, quindi:

$$\sum_{j=1}^{s} e_j \frac{(p_j-1)}{2}\sum_{i=1}^{r} f_i \frac{(q_i-1)}{2} = \left(\frac{(m-1)}{2} + 2k\right)\sum_{j=1}^{s} e_j \frac{(p_j-1)}{2} = \left(\frac{(m-1)}{2} + 2k\right)\left(\frac{(n-1)}{2} + 2h\right) =$$

$$= \frac{(m-1)}{2}\frac{(n-1)}{2} + 2\left(k\frac{(n-1)}{2} + h\frac{(m-1)}{2} + hk\right)$$

cioè:

$$\sum_{j=1}^{s} e_j \frac{(p_j-1)}{2}\sum_{i=1}^{r} f_i \frac{(q_i-1)}{2} \equiv \frac{(n-1)}{2}\frac{(m-1)}{2} \pmod 2$$ □

*Esempio 6.4-**Dirichlet fonte [ANT] Sec.9 Exercise 7***

Valutiamo $\left(\frac{365}{1847}\right)$ usando solo le proprietà del simbolo di Legendre, e usando il simbolo di Jacobi. Si noti che $1847$ è primo.

$$\left(\frac{365}{1847}\right) = \left(\frac{5}{1847}\right)\left(\frac{73}{1847}\right) = \left(\frac{1847}{5}\right)(-1)^{2\cdot 923}\left(\frac{1847}{73}\right)(-1)^{36\cdot 923} = \left(\frac{1847}{5}\right)\left(\frac{1847}{73}\right)$$

$$= \left(\frac{2+368\cdot 5}{5}\right)\left(\frac{22+25\cdot 73}{73}\right) = \left(\frac{2}{5}\right)\left(\frac{22}{73}\right) = (-1)^3\left(\frac{22}{73}\right) = -\left(\frac{22}{73}\right) = -\left(\frac{11}{73}\right)\left(\frac{2}{73}\right)$$

$$= -\left(\frac{73}{11}\right)(-1)^{5\cdot 36}(-1)^{666} = -\left(\frac{73}{11}\right) = -\left(\frac{7+11\cdot 6}{11}\right) = -\left(\frac{7}{11}\right) = -\left(\frac{11}{7}\right)(-1)^{5\cdot 3} = \left(\frac{11}{7}\right) =$$

$$= \left(\frac{4+7}{11}\right) = \left(\frac{4}{7}\right) = 1$$

Applichiamo ora il simbolo di Jacobi:

$$\left(\frac{365}{1847}\right) = (-1)^{182\cdot 923}\left(\frac{1847}{365}\right) = \left(\frac{22+5\cdot 365}{365}\right) = \left(\frac{22}{365}\right) = \left(\frac{2}{365}\right)\left(\frac{11}{365}\right) = (-1)^{16653}\left(\frac{11}{365}\right)$$

$$= -(-1)^{182\cdot 5}\left(\frac{365}{11}\right) = -\left(\frac{2+33*11}{11}\right) = -\left(\frac{2}{11}\right) = -(-1)^{15} = 1$$

## 7. *Residui quadratici modulo un numero intero.*

Si può ora affrontare per $n$ intero positivo la risoluzione della congruenza:

$$X^2 \equiv a \pmod{n} \qquad (7.1)$$

estendendo in modo naturale il concetto di residuo quadratico di un numero intero qualsiasi.

*Definizione 7.1.*
Sia $n$ un intero positivo ed $a$ un intero coprimo con $n$. Diciamo che $a$ è un residuo quadratico di $n$ quando la $X^2 \equiv a \pmod{n}$ ha almeno una soluzione, altrimenti $a$ è un residuo non quadratico.

Per gli sviluppi successivi conviene studiare le condizioni generali affinché sia risolvibile la congruenza:

$$f(X) \equiv 0 \pmod{n} \qquad (7.2)$$

dove $f(X)$ è un polinomio non nullo a coefficienti interi, il caso che a noi interessa corrisponderà a $f(X) = X^2 - a$.

*Proposizione 7.2*
Sia $n = p_1^{e_1} p_2^{e_2} \ldots p_s^{e_s}$ con $p_i$ primo, $e_i \geq 1$ ed $s \geq 1$. Le soluzioni della congruenza $f(X) \equiv 0 \pmod{n}$ coincidono con le soluzioni del sistema di congruenze:

$$\begin{cases} f(X) \equiv 0 \pmod{p_i^{e_i}} \\ 1 \leq i \leq s \end{cases} \qquad (7.3)$$

Dimostrazione
Se $x'$ è una soluzione di $f(X) \equiv 0 \pmod{n}$, dato che $p_i^{e_i} \mid n$ segue che $x'$ è soluzione di $f(X) \equiv 0 \pmod{p_i^{e_i}}$ $\forall i(1 \leq i \leq s)$.
Viceversa se $x''$ è una soluzione di $f(X) \equiv 0 \pmod{p_i^{e_i}}$ $\forall i(1 \leq i \leq s)$, e poiché $p_i^{e_i} \mid f(x'')$ $\forall i(1 \leq i \leq s)$ dato che $MCD(p_j^{e_j}, p_i^{e_i}) = 1$ per $j \neq i$, segue che $n \mid f(x'')$ e quindi $x''$ è una soluzione di $f(X) \equiv 0 \pmod{n}$. □

*Proposizione 7.3*
Sia $n = p_1^{e_1} p_2^{e_2} \ldots p_s^{e_s}$ con $p_i$ primo, $e_i \geq 1$ ed $s \geq 1$, se $\forall i(1 \leq i \leq s)$ $f(X) \equiv 0 \pmod{p_i^{e_i}}$ ammette $x_{ij}$ $1 \leq j \leq N_i$ soluzioni, allora la congruenza $f(X) \equiv 0 \pmod{n}$ ammette esattamente $N = \prod_{i=1}^{s} N_i$ soluzioni non congrue modulo $n$, ciascuna definita come segue:

$x_k = \sum_{i=1}^{s} x_{ij_{ik}} n_i \overline{n}_i$, $1 \leq j_{ik} \leq N_i$, $1 \leq k \leq N$, dove:

$n_i = \dfrac{n}{p_i^{e_i}}$, e $\overline{n}_i$ è un inverso moltiplicativo modulo $p_i^{e_i}$ di $n_i$, cioè $n_i \overline{n}_i \equiv 1 \pmod{p^{e_i}}$.

Dimostrazione
Se consideriamo una $s-upla$ di soluzioni $\left(x_{1j_{1k}}, x_{2j_{2k}}, \ldots, x_{sj_{sk}}\right)$, rispettivamente delle congruenze $f(X) \equiv 0 \pmod{p_1^{e_1}}$, $f(X) \equiv 0 \pmod{p^{e_2}}, \ldots, f(X) \equiv 0 \pmod{p^{e_s}}$ allora applicando il teorema cinese dei resti (vedi ad esempio [INT] §8.1) sappiamo che per il sistema:

$$\begin{cases} X \equiv x_{1j_1} \pmod{p_1^{e_1}} \\ X \equiv x_{2j_2} \pmod{p_2^{e_2}} \\ \quad .... \\ X \equiv x_{1j_s} \pmod{p_s^{e_s}} \end{cases}$$

esiste una sola classe di soluzioni tra loro congrue modulo $n$, è un suo rappresentante è $x_k = \sum_{i=1}^{s} x_{ij_{ik}} n_i \overline{n}_i$ con $n_i = \frac{n}{p_i^{e_i}}$, e $\overline{n}_i$ uno degli inversi moltiplicativi modulo $p_i^{e_i}$ di $n_i$, cioè $n_i\overline{n}_i \equiv 1 \pmod{p_i^{e_i}}$.

Siccome $\forall i(1 \le i \le s)$, $x_k \equiv x_{ij_{ik}} \pmod{p_1^{e_i}}$ quindi $f(x_k) \equiv f(x_{ij_{ik}}) \equiv 0 \pmod{p_i^{e_i}}$, segue che $x_k$ è una soluzione del sistema (7.3) e per la Proposizione 7.2 è anche una soluzione della congruenza (7.2). Percorrendo in senso inverso il nostro ragionamento da una soluzione $x_k$ di (7.2), si giunge ad una $s-upla$ di soluzioni $(x_{1j_{1k}}, x_{2j_{2k}}, \dots, x_{sj_{sk}})$ di (7.3). Per completare la dimostrazione basta notare infine che due soluzioni $x_k'$ e $x_k''$ congrue modulo $n$, danno luogo a due $s-uple$ le cui componenti sono congrue modulo $p_i^{e_i}$ $(1 \le i \le s)$, e viceversa. □

Ritorniamo ora ad occuparci del caso relativo alla congruenza (7.2) e consideriamo separatamente i due casi, $p = 2$ e $p$ dispari.

*Proposizione 7.3*

Sia $p$ un primo dispari ed $a$ un intero tale che $MCD(a, p) = 1$, allora la congruenza $X^2 \equiv a \pmod{p^e}$ con $e \ge 1$ è risolubile se e solo se $\left(\frac{a}{p}\right) = 1$ e ammette esattamente due soluzioni non congruenti modulo $p^e$.

Dimostrazione

Se la congruenza è risolubile allora lo è anche $X^2 \equiv a \pmod{p}$ e quindi $\left(\frac{a}{p}\right) = 1$.

Viceversa sia $\left(\frac{a}{p}\right) = 1$, procediamo per induzione su $e$. La congruenza è risolubile per ipotesi con $e = 1$ e ha esattamente due soluzioni non congruenti modulo $p$. Assumiamo che la congruenza sia risolubile per $e \ge 1$ e abbia esattamente due soluzioni $x_1$ e $x_2$ non congruenti modulo $p^e$, e consideriamo la congruenza $X^2 \equiv a \pmod{p^{e+1}}$. Sia $x_1^2 \equiv a \pmod{p^e}$ allora esiste $l_1 \in \mathbb{Z}$ tale che $x_1^2 = a + l_1 p^e$ inoltre $MCD(p, 2x_1) = 1$ quindi la congruenza $2x_1 Y \equiv -l_1 \pmod{p}$ ammette un'unica soluzione $y_1$.

Allora abbiamo che $x_1' = x_1 + y_1 p^e$ è una soluzione, infatti:

$x_1'^2 = x_1^2 + 2x_1 y_1 p^e + y_1^2 p^{2e} = a + l_1 p^e + (-l_1 + h_1 p)p^e + (y_1^2 p^{e-1})p^{e+1} \equiv a \pmod{p^{e+1}}$.

Analogamente ripetendo il ragionamento per $x_2$, abbiamo $x_2' = x_2 + y_2 p^e$ è una seconda soluzione e siccome se per assurdo fosse $p^{e+1} \mid (x_2' - x_1')$ allora in particolare $p^e \mid (x_2' - x_1')$ e quindi $p^e \mid (x_2 - x_1)$ da cui l'assurdo. Quindi $x_1'$ e $x_2'$ sono non congruenti modulo $p^{e+1}$. Sia ora $x'$ una soluzione di $X^2 \equiv a \pmod{p^{e+1}}$, allora $x' \equiv x_i' \pmod{p^e}$ con $i = 1$ oppure $i = 2$, sia quindi $x' = x_i + hp^e$, allora abbiamo che:

$x'^2 \equiv a \equiv x_i'^2 \pmod{p^{e+1}}$

$x_i^2 + 2x_i h p^e + h^2 p^{2e} \equiv x_i^2 + 2x_i y_i p^e + y_i^2 p^{2e} \pmod{p^{e+1}}$

$hp^e \equiv y_i p^e \pmod{p^{e+1}}$

$h \equiv y_i \pmod{p}$ da cui segue infine che $x' \equiv x_i' \pmod{p^{e+1}}$. □

*Proposizione 7.4*

La congruenza $X^2 \equiv a \pmod{2}$ è sempre risolubile per $a$ dispari e ha esattamente una sola classe di soluzioni congruenti modulo $2$, e cioè la classe dei numeri dispari.

Dimostrazione
Lasciata come esercizio al lettore.

*Proposizione 7.5*

La congruenza $X^2 \equiv a \pmod{4}$ con $MCD(a,4)=1$ è risolubile se e solo se $a \equiv 1 \pmod{4}$ e ha esattamente due soluzioni non congruenti modulo $4$.

Dimostrazione
Sia $x_1 \in \mathbb{Z}$ una soluzione di $X^2 \equiv a \pmod{4}$, dato che $a$ è dispari allora anche $x_1$ è dispari e dato che il quadrato di ogni intero e congruente a $1 \pmod{4}$, si ha che $x^2 \equiv a \equiv 1 \pmod{4}$.
Viceversa, se $a \equiv 1 \pmod{4}$ allora 1 e 3 sono le soluzioni della congruenza. □

*Proposizione 7.6*

Consideriamo la congruenza $X^2 \equiv a \pmod{2^e}$ con $e \geq 3$ e $MCD(a,8)=1$, valgono le seguenti.

**(i)** La congruenza è risolubile se e solo se $a \equiv 1 \pmod{8}$.

**(ii)** Se la congruenza è risolubile ha esattamente quattro soluzioni non congruenti modulo $2^e$, indicata con $x$ una di esse le altre sono $-x$, $x+2^{e-1}$, $-x+2^{e-1}$.

Dimostrazione
**Parte (i)**
Abbiamo che $(2k+1)^2 = 4k^2+4k+1 = 4k(k+1)+1 \equiv 1 \pmod{8}$, cioè il quadrato di ogni intero dispari e congruente ad $1 \pmod{8}$. Quindi dato che se $X^2 \equiv a \pmod{2^e}$ ammette soluzioni allora anche $X^2 \equiv a \pmod{8}$ ammette soluzioni, segue che $a \equiv 1 \pmod{8}$.
Viceversa sia $a \equiv 1 \pmod{8}$ procediamo per induzione su $e$. Se $e=3$ la congruenza è risolubile ed ha le soluzioni $1,3,5,7 \pmod{8}$. Supponiamo ora che per $e \geq 3$ la congruenza sia risolubile con quattro soluzioni $x_i$ con $i=1,2,3,4$ non congruenti modulo $8$, e dimostriamo che lo stesso vale per $e+1$. Fissiamo $i$, se $x_i$ è una soluzione di $X^2 \equiv a \pmod{2^e}$, allora esiste $l_i \in \mathbb{Z}$ tale che $x_i^2 = a + l_i 2^e$, inoltre $a$ e $x_i$ sono entrambi dispari per cui la congruenza lineare $x_i Y \equiv -l_i \pmod{2}$ ammette un'unica soluzione $y_i$ modulo $2$. Poniamo $x_i' = x_i + y_i 2^{e-1}$ allora risulta che $x_i'^2 \equiv a \pmod{2^{e+1}}$:

$$x_i'^2 = x_i^2 + 2x_i y_i 2^{e-1} + y_i^2 2^{2(e-1)} = a + l_i 2^e - l_i 2^e + u2^{e+1} + y_i^2 2^{2(e-1)} \equiv a \pmod{2^{e+1}}.$$

**Parte (ii)**
Potremmo ripetere lo stesso ragionamento fatto nella parte **(i)** per ogni $1 \leq i \leq 4$, però in generale non si ottengono quattro soluzioni non congruenti modulo $2^{e+1}$, infatti abbiamo:
$x_j' - x_i' = x_j - x_i + (y_j - y_i)2^{e-1}$ e risulta

$$2^{e+1} \mid \left(x_j' - x_i'\right) \text{ se e solo se } x_j - x_i \equiv 0 \pmod{2^{e-1}} \text{ e } 4 \mid \left(\frac{x_j - x_i}{2^{e-1}} + y_j - y_i\right).$$

Questo ad esempio succede segliendo opportunamente $y_j - y_i$ per $e=3$, considerando $x_j = 7$ e $x_i = 3$. Mentre, semplicemente svolgendo i calcoli, si verifica che se $x^2 \equiv a \pmod{2^e}$, allora sono soluzioni non congruenti modulo $2^e$: $x$, $-x$, $x+2^{e-1}$, $-x+2^{e-1}$.

Quindi in particolare scelto $i$ con $1 \le i \le 4$ procedendo come definito nella parte **(i)**, sono soluzioni non congruenti modulo $2^{e+1}$ della $X^2 \equiv a \pmod{2^{e+1}}$ le seguenti: $x_i'$, $x_i' + 2^e$, $-x_i'$, $-x_i' + 2^e$.
Procediamo di nuovo per induzione, per $e = 3$ i residui quadratici sono esattamente $4$ per cui il risultato è immediato. Supponiamo ora **(ii)** vera per $e$ e dimostriamola per $e+1$.
Sia $x'$ tale che $x'^2 \equiv a \pmod{2^{e+1}}$ è sempre possibile ricondursi al caso $1 \le x' < 2^{e+1}$, in particolare $x'^2 \equiv a \pmod{2^e}$ quindi per qualche $1 \le i \le 4$ risulta $x' \equiv x_i \pmod{2^e}$, e di nuovo non è restrittivo supporre che $1 \le x_i < 2^e$, segue che si hanno due possibilità:

- $x' = x_i$, e quindi $x_i^2 \equiv a \pmod{2^{e+1}}$, quindi $x_i' = x_i$ e $x' \equiv x_i' \pmod{2^{e+1}}$.
  Quindi $-x' \equiv -x_i' \pmod{2^{e+1}}$, $\pm x' + 2^e \equiv \pm x_i' + 2^e \pmod{2^{e+1}}$
- $x' = x_i + 2^e$, allora:
  - se $x_i' = x_i$ ne segue che $x' = x_i' + 2^e$
  - se invece $x_i' = x_i + 2^{e-1}$ allora essendo $x'^2 \equiv a \equiv x_i'^2 \pmod{2^{e+1}}$ segue che si ha l'assurdo $x_i \equiv 0 \pmod 2$.

  Quindi $-x' \equiv -x_i' - 2^e = -x_i' + 2^e - 2^{e+1} \equiv -x_i' + 2^e \pmod{2^{e+1}}$ e
  $\pm x' + 2^e \equiv \pm x_i' + 2^{e+1} \equiv \pm x_i' \pmod{2^{e+1}}$. □

Siamo a questo punto pronti ad enunciare il teorema generale relativo alla condizioni generali affinché la congruenza $X^2 \equiv a \pmod n$ sia risolvibile.

*Teorema 7.7*

Sia $n = 2^{e_0} p_1^{e_1} \dots p_s^{e_s}, a \in \mathbb{Z}$ con $p_1, \dots, p_s$ numeri primi dispari distinti e $MCD(a,n) = 1$, allora la congruenza $X^2 \equiv a \pmod n$:

**(i)** è risolubile se e solo se:

- $\left(\dfrac{a}{p_1}\right) = \dots = \left(\dfrac{a}{p_s}\right) = 1$
- $\begin{cases} a \equiv 1 \pmod 2 & \text{se } e_0 = 1 \\ a \equiv 1 \pmod 4 & \text{se } e_0 = 2 \\ a \equiv 1 \pmod 8 & \text{se } e_0 \ge 3 \end{cases}$

**(ii)** Se la congruenza è risolubile il numero delle sue soluzioni incongruenti modulo $n$ è dato da:

$$\begin{cases} 2^s & \text{se } e_0 \le 1 \\ 2^{s+1} & \text{se } e_0 = 2 \\ 2^{s+2} & \text{se } e_0 \ge 3 \end{cases}$$

Dimostrazione
Diretta conseguenza delle Proposizioni da 7.2 a 7.6. □

Legendre, come Gauss stesso cita nell'articolo 151 di [DA], era giunto ad un teorema che è equivalente al *"Theorema Fundamentale"*

*Proposizione 7.8*-***Formulazione di Legendre della reciprocità quadratica***

Se $p$ e $q$ sono due numeri primi dispari, allora:

i. Se $p \equiv 1 \pmod 4$ o $q \equiv 1 \pmod 4$, allora $p^{\frac{q-1}{2}} \equiv \varepsilon \pmod q$, $q^{\frac{p-1}{2}} \equiv \varepsilon \pmod p$, con $\varepsilon \in \{-1,+1\}$.

ii. Se $p \equiv 3 \pmod 4$ e $q \equiv 3 \pmod 4$, allora $p^{\frac{q-1}{2}} \equiv \varepsilon \pmod q$, $q^{\frac{p-1}{2}} \equiv -\varepsilon \pmod p$, con $\varepsilon \in \{-1,+1\}$.

Dimostrazione di Legendre
Sempre citando [DA] articolo 151 a proposito della dimostrazione di Legendre della legge della reciprocità quadratica*"…egli presupponeva molte cose senza dimostrazione (…), alcune delle quali non sono state dimostrate da nessuno fino ad adesso, e alcune delle quali a nostro giudizio non sono dimostrabili senza l'aiuto del teorema fondamentale,….,e quindi la nostra dimostrazione deve essere considerata come la prima…."*

Vediamo con maggiore dettaglio le "cose" presupposte. Legendre ha dimostrato la seguente proposizione:

*"Condizione necessaria e sufficiente affinché l'equazione diofantea* $fX^2 + gY^2 = hZ^2$ *con* $f,g,h \in \mathbb{Z}$*, coprimi a coppie e privi di divisori quadratici, sia risolvibile e che siano risolvibili le tre congruenze quadratiche:*
$T^2 \equiv -fg \pmod h$, $T^2 \equiv hg \pmod f$, $T^2 \equiv sf \pmod g$ *"*

La verifica del fatto che la condizione è necessaria e relativamente semplice, infatti ad esempio per ogni primo $p \mid h$:

$$\left(\frac{-fg}{p}\right) = \left(\frac{-fgx_0^2}{p}\right) = \left(\frac{g^2y_0^2 - ghz_0^2}{p}\right) = \left(\frac{g^2y_0^2}{p}\right) = 1$$

che per il teorema 7.7 equivale a che sia risolvibile $T^2 \equiv -fg \pmod h$, analogamente si giunge alle altre due congruenze, per la dimostrazione di sufficienza si veda [DE] capitolo 7 teorema 3, oppure [DA] articolo 294.

Con l'utilizzo della proposizione sopra citata e facendo alcune assunzioni arbitrarie, anche se come vedremo intuitive, Legendre è riuscito a dimostrare la legge della reciprocità quadratica considerando separatamente i seguenti cinque casi:

*Caso I) Se* $p$ *e* $q$ *sono primi dispari distinti della forma* $4n+3$*, allora non si può avere* $\left(\frac{p}{q}\right) = 1$ *e* $\left(\frac{q}{p}\right) = 1$.

Infatti assumendo $f = p$, $g = q$ e $h = 1$ abbiamo che in base alla proposizione sopra enunciata $pX^2 + qY^2 = Z^2$ ammette soluzioni non banali. Una generica soluzione $(x_0, y_0, z_0)$ è sempre riconducibile ad una soluzione costituita da coppie di numeri coprimi a coppie, in cui al più uno dei tre numeri può essere pari, quindi $px_0^2 + qy_0^2 - z_0^2 \equiv 2 \pmod 4$ da cui l'assurdo.

*Caso II) se* $p$ *è un numero primo della forma* $4n+1$*, e* $q$ *è un numero primo della forma* $4n+3$ *allora non possiamo avere* $\left(\frac{p}{q}\right) = -1$ *e* $\left(\frac{q}{p}\right) = 1$

Infatti altrimenti si avrebbe $\left(\frac{-p}{q}\right)=1$ e quindi l'equazione $X^2 + pY^2 = qZ^2$ ammetterebbe soluzioni non banali, ma essendo $x_0{}^2 + py_0{}^2 - qz_0{}^2 \equiv 2 \pmod 4$ questo è impossibile.

*Caso III) Se $p$ e $q$ sono primi dispari distinti della forma $4n+1$, allora non si può avere* $\left(\frac{p}{q}\right)=1$ *e* $\left(\frac{q}{p}\right)=-1$.

Se assumiamo l'esistenza di un altro numero primo $r$ della forma $4n+3$ tale che $\left(\frac{r}{q}\right)=1$ e $\left(\frac{p}{r}\right)=-1$, allora mediante (II) abbiamo che $\left(\frac{q}{r}\right)=1$, $\left(\frac{r}{p}\right)=-1$. Se ora per assurdo assumiamo $\left(\frac{p}{q}\right)=1$ e $\left(\frac{q}{p}\right)=-1$, allora abbiamo che $\left(\frac{qr}{p}\right)=1$, $\left(\frac{pr}{q}\right)=1$, $\left(\frac{pq}{r}\right)=-1$ da cui essendo $r$ della forma $4n+3$ segue che $\left(\frac{-pq}{r}\right)=1$. Quindi $pX^2 + qY^2 = rZ^2$ è risolvibile e questo è assurdo dato che $px_0{}^2 + qy_0{}^2 - rz_0{}^2 \equiv 2 \pmod 4$.

*Caso IV) se $p$ è un numero primo della forma $4n+1$, e $q$ è un numero primo della forma $4n+3$ allora non possiamo avere* $\left(\frac{p}{q}\right)=1$ *e* $\left(\frac{q}{p}\right)=-1$.

Assumiamo l'esistenza di un numero primo ausiliario $r$ della forma $4n+1$, tale che $\left(\frac{r}{p}\right)=\left(\frac{r}{q}\right)=-1$. Allora abbiamo tramite (II) che $\left(\frac{q}{r}\right)=-1$ e tramite (III) che $\left(\frac{p}{r}\right)=-1$, quindi $\left(\frac{pq}{r}\right)=1$. Se per assurdo fosse $\left(\frac{p}{q}\right)=1$ e $\left(\frac{q}{p}\right)=-1$, allora si avrebbe $\left(\frac{pr}{q}\right)=-1$, $\left(\frac{-pr}{q}\right)=1$ e $\left(\frac{qr}{p}\right)=1$, quindi l'equazione $pX^2 + rY^2 = qZ^2$ sarebbe risolvibile in contraddizione con il fatto che $px_0{}^2 + ry_0{}^2 - qz_0{}^2 \equiv 2 \pmod 4$.

*Caso V) Se $p$ e $q$ sono primi dispari distinti della forma $4n+3$, allora non si può avere* $\left(\frac{p}{q}\right)=-1$ *e* $\left(\frac{q}{p}\right)=-1$.

Sia per assurdo $\left(\frac{p}{q}\right)=-1$ e $\left(\frac{q}{p}\right)=-1$, assumiamo l'esistenza di un numero primo ausiliario $r$ della forma $4n+1$, tale che $\left(\frac{r}{p}\right)=\left(\frac{r}{q}\right)=-1$. Allora abbiamo che $\left(\frac{qr}{p}\right)=\left(\frac{pr}{q}\right)=1$, inoltre tramite (II) $\left(\frac{p}{r}\right)=\left(\frac{q}{r}\right)=-1$ e quindi $\left(\frac{pq}{r}\right)=\left(\frac{-pq}{r}\right)=1$, pertanto l'equazione $pX^2 + qY^2 = rZ^2$ è risolvibile in contraddizione con $px_0{}^2 + qy_0{}^2 - rz_0{}^2 \equiv 2 \pmod 4$.

L'assumere nei casi III-IV-V l'esistenza di numeri ausiliari che soddisfano particolari condizioni relativamente a $p$ e $q$ rende la dimostrazione per questa via incompleta. □

*Osservazione 7.9*

Siamo ora in grado di approfondire quanto detto nel paragrafo 5, si consideri $s$ primo dispari e il primo $q$, tale che $\left(\frac{q}{s}\right)=1$, allora almeno una delle due equazioni diofantee

$$qX^2 - Y^2 = sZ^2 \tag{7.3}$$

$$X^2 - qY^2 = sZ^2 \tag{7.4}$$

ammette soluzioni. Infatti per entrambe è verificata $T^2 \equiv q \pmod{s}$, mentre deve risultare $\left(\frac{-s}{q}\right)=1$ affinché la (7.3) sia risolvibile, e $\left(\frac{s}{q}\right)=1$ affinché la (7.4) sia risolvibile.

Se $\left(\frac{s}{q}\right)=\left(\frac{-s}{q}\right)$, allora $q \equiv 1 \pmod{4}$ e allora tenuto conto che $\left(\frac{q}{s}\right)=1$ per la reciprocità quadratica $\left(\frac{s}{q}\right)=\left(\frac{-s}{q}\right)=1$, pertanto almeno una tra le (7.3) e le (7.4) è risolvibile.

*Osservazione 7.10*

Può succedere che $a$ non è un residuo quadratico modulo $n$, però $\left(\frac{a}{n}\right)=1$. Ad esempio $\left(\frac{2}{9}\right)=\left(\frac{2}{3}\right)^2=1$ ma $2$ non è un quadrato modulo $9$.

Viceversa se $a$ è un residuo quadratico modulo $n$, allora $\left(\frac{a}{n}\right)=1$ come si deduce immediatamente dalla definizione del simbolo di Jacobi e dal Teorema 7.7.

## 8. *La congruenza quadratica* $aX^2 + bX + c \equiv 0 \pmod{n}$

Siamo ora in grado di analizzare:

$$aX^2 + bX + c \equiv 0 \pmod{n} \qquad (8.1)$$

Ovviamente nel caso $a \equiv 0 \pmod{n}$ si ricade nel caso delle congruenze lineari non oggetto di questa nota (vedi [INT]). Moltiplicando per $4a$ ambo i membri della (8.1) otteniamo la congruenza:

$$4a^2X^2 + 4abX + 4ac \equiv 0 \pmod{4an} \qquad (8.2)$$

La (8.2) è equivalente alla (8.1) cioè ogni numero che soddisfa una delle due soddisfa anche l'altra. Infatti per ogni $x \in \mathbb{Z}$ $4an \mid \left(4a^2x^2 + 4abx + 4ac\right)$ equivale $n \mid \left(ax^2 + bx + c\right)$.

La congruenza (8.2) può essere riscritta come:

$$(2aX + b)^2 \equiv (b^2 - 4ac) \pmod{4an} \qquad (8.3)$$

ovvero ci siamo ricondotti alla risoluzione della congruenza

$$T^2 \equiv b^2 - 4ac \pmod{4an} \qquad (8.4)$$

considerando di questa una soluzioni $t$ dobbiamo risolvere la congruenza:

$$2aX + b \equiv t \pmod{4an} \qquad (8.5)$$

questa a sua volta è risolubile (vedi [INT] Teorema 57) se e solo se $2a \mid (t - b)$ infatti $MCD(2a, 4an) = 2a$, quindi se $t \equiv b \pmod{2a}$.

*Teorema 8.1*

Condizione necessaria affinché la congruenza (8.1) sia risolubile e che il discriminante $\Delta = b^2 - 4ac$ sia un residuo quadratico modulo $4an$. Ad una soluzione $t$ della congruenza $T^2 \equiv \Delta \pmod{4ac}$ corrisponde una soluzione della congruenza (8.1) se e solo se $t \equiv b \pmod{2a}$.

Dimostrazione

Segue dai ragionamenti svolti preliminarmente all'enunciazione di questo teorema. □

Nel caso particolare di $MCD(2a, n) = 1$ è possibile ottenere un risultato più generale, espresso dal seguente teorema.

*Teorema 8.2*

**(i)** Se $MCD(2a, n) = 1$ allora condizione necessaria e sufficiente affinché la congruenza (8.1) sia risolubile e che il discriminante $\Delta = b^2 - 4ac$ sia un residuo quadratico di modulo $n$.

**(ii)** Le soluzioni della (8.1) non congruenti modulo $n$ sono in numero uguale alle soluzioni della congruenza

$$T^2 \equiv \Delta \pmod{n}. \qquad (8.6)$$

**(iii)** Il legame tra una soluzione della (8.1) ed una soluzione della (8.6) è dato da:

$$x = \frac{n+1}{2}\bar{a}(t - b) \qquad (8.8)$$

dove $\bar{a}$ è un inverso moltiplicativo modulo $n$ di $a$.

Dimostrazione

Che la condizione è necessaria segue dal teorema 8.1 dato che se $T^2 \equiv \Delta \pmod{4ac}$ è risolubile allora anche $T^2 \equiv \Delta \pmod{n}$ è risolubile.

Viceversa se $T^2 \equiv \Delta \pmod{n}$ è risolubile allora sia $t$ una sua soluzione, consideriamo la congruenza lineare:

$2aX + b \equiv t \pmod{n}$ siccome $MCD(2a, n) = 1$ questa ammette una ed una sola soluzione $x$ modulo $n$, dato che:

$4a(ax^2 + bx + c) = (2ax + b)^2 - (b^2 - 4ac) \equiv 0 \pmod{n}$, cioè $n \mid 4a(ax^2 + bx + c)$ e inoltre $MCD(2a, n) = 1$, segue che $n \mid (ax^2 + bx + c)$ cioè $ax^2 + bx + c \equiv 0 \pmod{n}$.

Da $MCD(2a,n)=1$ segue che $MCD(a,n)=1$ e $MCD(2,n)=1$ (cioè $n$ è dispari), quindi $a$ ha un inverso moltiplicativo modulo $n$ cioè un numero $\overline{a}$ tale che $a\overline{a}\equiv 1 \pmod n$. Analogamente 2 ammette un inverso moltiplicativo infatti esso è $\frac{n+1}{2}$.

Pertanto se esiste $x$ tale che: $2ax+b\equiv t \pmod n$ allora $x\equiv \frac{n+1}{2}\overline{a}2ax=\frac{n+1}{2}\overline{a}(t-b) \pmod n$.

Infine se abbiamo due soluzioni $t'$ e $t''$ di (8.6),allora:

$x''-x'\equiv \frac{n+1}{2}\overline{a}(t''-t') \pmod n$ siccome $2a\frac{n+1}{2}\overline{a}=1+kn$ e $MCD(2a,n)=1$ segue che $n \mid (t''-t')$ da cui $t''\equiv t' \pmod n$. □

*Esempio 8.3*

Come esempio risolviamo con in metodi del presente paragrafo e del paragrafo precedente, la seguente congruenza:

$3X^2+7X-1\equiv 0 \pmod{15}$

chiaramente in questo caso per tentativi sugli interi da 1 a 14 si arriverebbe più rapidamente alla soluzione, ma lo scopo di questo esempio è quello di dare una visione generale di quanto esposto.

Dobbiamo risolvere la congruenza quadratica:

$T^2\equiv 61 \pmod{2^2\cdot 3^2\cdot 5}$

risulta che:

$61\equiv 1 \pmod 4$

$\left(\frac{61}{3}\right)=\left(\frac{1+20\cdot 3}{3}\right)=1$

$\left(\frac{61}{5}\right)=\left(\frac{1+12\cdot 5}{3}\right)=1$

A questo punto risolviamo le tre congruenze per 4, 9 e 5.

$T^2\equiv 61 \pmod 4$

$t_{11}=1$, $t_{12}=3$.

$T^2\equiv 61 \pmod 9$ equivalente a $T^2\equiv 7 \pmod 9$

risolviamo prima $T^2\equiv 7 \pmod 3$ per cui abbiamo come soluzioni $t'_{21}=1$, $t'_{22}=2$, quindi

${t'_{21}}^2=7-2\cdot 3$ risolviamo al congruenza $2Y\equiv 2 \pmod 3$ con unica soluzione $y_{21}=1$.

$t'_{22}=7-1\cdot 3$ risolviamo al congruenza $4Y\equiv 1 \pmod 3$ con unica soluzione $y_{21}=1$.

Segue che le soluzioni sono $t_{21}=1+1\cdot 3=4$, $t_{22}=2+1\cdot 3=5$

$T^2\equiv 61 \pmod 5$ equivalente a $T^2\equiv 1 \pmod 5$ per cui le soluzioni sono $t_{31}=1$, $t_{32}=4$

Usando le notazioni della proposizione xx abbiamo la tabella:

| $p_i^{e_i}$ | $t_{i1}$ | $t_{i2}$ | $n_i$ | $\overline{n}_i$ |
|---|---|---|---|---|
| **4** | 1 | 3 | 45 | 1 |
| **9** | 4 | 5 | 20 | 5 |
| **5** | 1 | 4 | 36 | 1 |

Segue si hanno le soluzioni:

$t_1=1\cdot 45\cdot 1+4\cdot 20\cdot 5+1\cdot 36\cdot 1\equiv 121 \pmod{180}$    $6\mid(121-7)$

$t_2=1\cdot 45\cdot 1+5\cdot 20\cdot 5+1\cdot 36\cdot 1\equiv 41 \pmod{180}$    $6\nmid(41-7)$

$t_3=1\cdot 45\cdot 1+5\cdot 20\cdot 5+4\cdot 36\cdot 1\equiv 149 \pmod{180}$    $6\nmid(149-7)$

$t_4=1\cdot 45\cdot 1+4\cdot 20\cdot 5+4\cdot 36\cdot 1\equiv 49 \pmod{180}$    $6\mid(49-7)$

$t_5=3\cdot 45\cdot 1+4\cdot 20\cdot 5+1\cdot 36\cdot 1\equiv 31 \pmod{180}$    $6\mid(31-7)$

$t_6 = 3 \cdot 45 \cdot 1 + 5 \cdot 20 \cdot 5 + 1 \cdot 36 \cdot 1 \equiv 131 \pmod{180} \quad 6 \nmid (131-7)$
$t_7 = 3 \cdot 45 \cdot 1 + 5 \cdot 20 \cdot 5 + 4 \cdot 36 \cdot 1 \equiv 59 \pmod{180} \quad 6 \nmid (59-7)$
$t_8 = 3 \cdot 45 \cdot 1 + 4 \cdot 20 \cdot 5 + 4 \cdot 36 \cdot 1 \equiv 139 \pmod{180} \quad 6 \mid (139-7)$

Si tratta ora di risolvere quattro congruenze lineari:
$6X \equiv 114 \pmod{180}$ ovvero $X \equiv 19 \pmod{30}$ da cui $x \equiv 4 \pmod{15}$
$6X \equiv 42 \pmod{180}$ ovvero $X \equiv 7 \pmod{30}$ e anche $x \equiv 7 \pmod{15}$
$6X \equiv 24 \pmod{180}$ ovvero $X \equiv 4 \pmod{30}$ e anche $x \equiv 4 \pmod{15}$
$6X \equiv 132 \pmod{180}$ ovvero $X \equiv 22 \pmod{30}$ da cui $x \equiv 7 \pmod{15}$
Si hanno quindi due soluzioni $4$ e $7$.

*Esempio 8.4*

$3X^2 + 7X - 1 \equiv 0 \pmod{3 \cdot 5 \cdot 13}$
è una semplice variazione dell'esempio precedente ma in questo caso la verifica a mano è meno economica. Utilizzando i calcoli già svolti nell'esempio precedente abbiamo che:

$$\left(\frac{61}{13}\right) = \left(\frac{9}{13}\right) = 1$$

risolviamo $T^2 \equiv 61 \pmod{13}$ che ha come soluzioni $3$ e $10$, abbiamo quindi che ogni soluzione della congruenza $T^2 \equiv 61 \pmod{2^2 \cdot 3^2 \cdot 5}$ a due soluzioni di $T^2 \equiv 61 \pmod{2^2 \cdot 3^2 \cdot 5 \cdot 13}$.

| $p_i^{e_i}$ | $t_{i1}$ | $t_{i2}$ | $n_i$ | $\overline{n}_i$ |
|---|---|---|---|---|
| **4** | 1 | 3 | 45*13 | 1*1 |
| **9** | 4 | 5 | 20*13 | 5*7 |
| **5** | 1 | 4 | 36*13 | 1*2 |
| **13** | 3 | 10 | 180 | 6 |

Quindi:
$t_1' = 1 \cdot 45 \cdot 13 \cdot 1 \cdot 1 + 4 \cdot 20 \cdot 13 \cdot 5 \cdot 7 + 1 \cdot 36 \cdot 13 \cdot 1 \cdot 2 + 3 \cdot 180 \cdot 6 \equiv 1381 \pmod{2340}$ e $6 \mid (1381-7)$
$t_1'' = 1 \cdot 45 \cdot 13 \cdot 1 \cdot 1 + 4 \cdot 20 \cdot 13 \cdot 5 \cdot 7 + 1 \cdot 36 \cdot 13 \cdot 1 \cdot 2 + 10 \cdot 180 \cdot 6 \equiv 1921 \pmod{2340}$ e $6 \mid (1921-7)$
$t_2' = 1 \cdot 45 \cdot 13 \cdot 1 \cdot 1 + 5 \cdot 20 \cdot 13 \cdot 5 \cdot 7 + 1 \cdot 36 \cdot 13 \cdot 1 \cdot 2 + 3 \cdot 180 \cdot 6 \equiv 1121 \pmod{2340}$ e $6 \nmid (1121-7)$
$t_2'' = 1 \cdot 45 \cdot 13 \cdot 1 \cdot 1 + 5 \cdot 20 \cdot 13 \cdot 5 \cdot 7 + 1 \cdot 36 \cdot 13 \cdot 1 \cdot 2 + 10 \cdot 180 \cdot 6 \equiv 1661 \pmod{2340}$ e $6 \nmid (1661-7)$
$t_3' = 1 \cdot 45 \cdot 13 \cdot 1 \cdot 1 + 5 \cdot 20 \cdot 13 \cdot 5 \cdot 7 + 4 \cdot 36 \cdot 13 \cdot 1 \cdot 2 + 3 \cdot 180 \cdot 6 \equiv 1589 \pmod{2340}$ e $6 \nmid (1589-7)$
$t_3'' = 1 \cdot 45 \cdot 13 \cdot 1 \cdot 1 + 5 \cdot 20 \cdot 13 \cdot 5 \cdot 7 + 4 \cdot 36 \cdot 13 \cdot 1 \cdot 2 + 10 \cdot 180 \cdot 6 \equiv 2129 \pmod{2340}$ e $6 \nmid (2129-7)$
$t_4' = 1 \cdot 45 \cdot 13 \cdot 1 \cdot 1 + 4 \cdot 20 \cdot 13 \cdot 5 \cdot 7 + 4 \cdot 36 \cdot 13 \cdot 1 \cdot 2 + 3 \cdot 180 \cdot 6 \equiv 1849 \pmod{2340}$ e $6 \mid (1849-7)$
$t_4'' = 1 \cdot 45 \cdot 13 \cdot 1 \cdot 1 + 4 \cdot 20 \cdot 13 \cdot 5 \cdot 7 + 4 \cdot 36 \cdot 13 \cdot 1 \cdot 2 + 10 \cdot 180 \cdot 6 \equiv 49 \pmod{2340}$ e $6 \mid (49-7)$
$t_5' = 3 \cdot 45 \cdot 13 \cdot 1 \cdot 1 + 4 \cdot 20 \cdot 13 \cdot 5 \cdot 7 + 1 \cdot 36 \cdot 13 \cdot 1 \cdot 2 + 3 \cdot 180 \cdot 6 \equiv 211 \pmod{2340}$ e $6 \mid (211-7)$
$t_5'' = 3 \cdot 45 \cdot 13 \cdot 1 \cdot 1 + 4 \cdot 20 \cdot 13 \cdot 5 \cdot 7 + 1 \cdot 36 \cdot 13 \cdot 1 \cdot 2 + 10 \cdot 180 \cdot 6 \equiv 751 \pmod{2340}$ e $6 \mid (751-7)$
$t_6' = 3 \cdot 45 \cdot 13 \cdot 1 \cdot 1 + 5 \cdot 20 \cdot 13 \cdot 5 \cdot 7 + 1 \cdot 36 \cdot 13 \cdot 1 \cdot 2 + 3 \cdot 180 \cdot 6 \equiv 2291 \pmod{2340}$ e $6 \nmid (2291-7)$
$t_6'' = 3 \cdot 45 \cdot 13 \cdot 1 \cdot 1 + 5 \cdot 20 \cdot 13 \cdot 5 \cdot 7 + 1 \cdot 36 \cdot 13 \cdot 1 \cdot 2 + 10 \cdot 180 \cdot 6 \equiv 491 \pmod{2340}$ e $6 \nmid (491-7)$
$t_7' = 3 \cdot 45 \cdot 13 \cdot 1 \cdot 1 + 5 \cdot 20 \cdot 13 \cdot 5 \cdot 7 + 4 \cdot 36 \cdot 13 \cdot 1 \cdot 2 + 3 \cdot 180 \cdot 6 \equiv 419 \pmod{2340}$ e $6 \nmid (419-7)$
$t_7'' = 3 \cdot 45 \cdot 13 \cdot 1 \cdot 1 + 5 \cdot 20 \cdot 13 \cdot 5 \cdot 7 + 4 \cdot 36 \cdot 13 \cdot 1 \cdot 2 + 10 \cdot 180 \cdot 6 \equiv 959 \pmod{2340}$ e $6 \nmid (959-7)$
$t'_8 = 3 \cdot 45 \cdot 13 \cdot 1 \cdot 1 + 4 \cdot 20 \cdot 13 \cdot 5 \cdot 7 + 4 \cdot 36 \cdot 13 \cdot 1 \cdot 2 + 3 \cdot 180 \cdot 6 \equiv 679 \pmod{2340}$ e $6 \mid (679-7)$
$t_8'' = 3 \cdot 45 \cdot 13 \cdot 1 \cdot 1 + 4 \cdot 20 \cdot 13 \cdot 5 \cdot 7 + 4 \cdot 36 \cdot 13 \cdot 1 \cdot 2 + 10 \cdot 180 \cdot 6 \equiv 1219 \pmod{2340}$ e $6 \mid (1219-7)$

Abbiamo infine le quattro soluzioni distinte:
$x_1 = 7$, $x_2 = 34$, $x_3 = 112$, $x_4 = 124$.

*Esempio 8.5*

$3X^2 + 7X - 1 \equiv 0 \pmod{5 \cdot 13 \cdot 19}$

è di nuovo una piccola variazione dell'equazione precedente dove ora $MCD(2a, n) = 1$. Dato che: $\left(\frac{61}{19}\right) = \left(\frac{4}{19}\right) = 1$ quindi $T^2 \equiv 61 \pmod{19}$ ha le soluzioni $2$ e $17$, quindi abbiamo la tabella riassuntiva:

| $p_i^{e_i}$ | $t_{i1}$ | $t_{i2}$ | $n_i$ | $\overline{n}_i$ |
|---|---|---|---|---|
| **5** | 1 | 4 | 247 | 3 |
| **13** | 3 | 10 | 95 | 10 |
| **19** | 2 | 17 | 65 | 12 |

La congruenza $T^2 \equiv 61 \pmod{5 \cdot 13 \cdot 19}$ ha quindi le otto soluzioni:

$t_1 = 1 \cdot 247 \cdot 3 + 3 \cdot 95 \cdot 10 + 2 \cdot 65 \cdot 12 \equiv 211 \pmod{5 \cdot 13 \cdot 19}$

$t_2 = 1 \cdot 247 \cdot 3 + 10 \cdot 95 \cdot 10 + 2 \cdot 65 \cdot 12 \equiv 686 \pmod{5 \cdot 13 \cdot 19}$

$t_3 = 1 \cdot 247 \cdot 3 + 10 \cdot 95 \cdot 10 + 17 \cdot 65 \cdot 12 \equiv 36 \pmod{5 \cdot 13 \cdot 19}$

$t_4 = 1 \cdot 247 \cdot 3 + 3 \cdot 95 \cdot 10 + 17 \cdot 65 \cdot 12 \equiv 796 \pmod{5 \cdot 13 \cdot 19}$

$t_5 = 4 \cdot 247 \cdot 3 + 3 \cdot 95 \cdot 10 + 2 \cdot 65 \cdot 12 \equiv 1199 \pmod{5 \cdot 13 \cdot 19}$

$t_6 = 4 \cdot 247 \cdot 3 + 10 \cdot 95 \cdot 10 + 2 \cdot 65 \cdot 12 \equiv 439 \pmod{5 \cdot 13 \cdot 19}$

$t_7 = 4 \cdot 247 \cdot 3 + 10 \cdot 95 \cdot 10 + 17 \cdot 65 \cdot 12 \equiv 1024 \pmod{5 \cdot 13 \cdot 19}$

$t_8 = 4 \cdot 247 \cdot 3 + 3 \cdot 95 \cdot 10 + 17 \cdot 65 \cdot 12 \equiv 549 \pmod{5 \cdot 13 \cdot 19}$

L'inverso moltiplicativo di $3$ modulo $1235$ è $412$, quindi abbiamo le seguenti otto soluzioni:

$x_1 = 412 \cdot 618 \cdot (211 - 7) \equiv 34 \pmod{1235}$

$x_2 = 412 \cdot 618 \cdot (686 - 7) \equiv 319 \pmod{1235}$

$x_3 = 412 \cdot 618 \cdot (36 - 7) \equiv 1034 \pmod{1235}$

$x_4 = 412 \cdot 618 \cdot (796 - 7) \equiv 749 \pmod{1235}$

$x_5 = 412 \cdot 618 \cdot (1199 - 7) \equiv 1022 \pmod{1235}$

$x_6 = 412 \cdot 618 \cdot (439 - 7) \equiv 72 \pmod{1235}$

$x_7 = 412 \cdot 618 \cdot (1024 - 7) \equiv 787 \pmod{1235}$

$x_8 = 412 \cdot 618 \cdot (549 - 7) \equiv 502 \pmod{1235}$

## 9. *Interi somma di due quadrati*

Ci occupiamo per $n \in \mathbb{Z}$ con $n > 0$, della soluzione dell'equazione diofantea:

$$X^2 + Y^2 = n \qquad (9.1)$$

Un primo risultato sull'argomento è dovuto a Fermat con la seguente proposizione:

*Proposizione 9.1*

Un numero primo $p$ è la somma di sue quadrati se e solo se $p = 2$ oppure $p \equiv 1 \pmod 4$.

Dimostrazione

Se $p = a^2 + b^2$ allora $a$ e $b$ non possono essere entrambi pari altrimenti $4 \mid p$. Se $a$ e $b$ sono entrambi dispari, allora $p \equiv 1 + 1 \pmod 4$, quindi $2 \mid p$ da cui $p = 2$. Se ad esempio $a$ è dispari e $b$ è pari, allora $p \equiv 1 + 0 = 1 \pmod 4$.

Viceversa, $2^2 = 1^2 + 1^2$, sia invece $p \equiv 1 \pmod 4$ allora $-1$ è un residuo quadratico modulo $p$ (Proposizione 3.28(c)), quindi esiste $x$, $1 \le x \le p-1$, tale che $x^2 + 1 = mp$, con $1 \le m \le p-1$. Quindi l'insieme $\{m \in \mathbb{Z} : 1 \le m \le p-1 \wedge \exists x, y \in \mathbb{Z}(mp = x^2 + y^2)\}$ non è vuoto. Sia $m_0$ il minimo intero in questo insieme, allora $1 \le m_0 \le p-1$. Dimostriamo che $m_0 = 1$ e come conseguenza $p$ è la somma di due quadrati. Assumiamo per assurdo che $1 < m_0$. Allora posto $x^2 + y^2 = m_0 p$ osserviamo scriviamo usando l'algoritmo della divisione:

$$\begin{cases} x = cm_0 + x_1 \\ y = dm_0 + y_1 \end{cases}$$

con $-m_0/2 < x_1, y_1 < m_0/2$ e $c$ e $d$ interi, inoltre $x_1 \ne 0 \vee y_1 \ne 0$, infatti $m_0 \nmid x \vee m_0 \nmid y$ altrimenti $m_0 = p$ che è assurdo. Abbiamo quindi che:

$0 < x_1^2 + y_1^2 < m_0^2/4 + m_0^2/4 = m_0^2/2 < m_0^2$ e $x_1^2 + y_1^2 \equiv x^2 + y^2 \equiv 0 \pmod{m_0}$ cioè $x_1^2 + y_1^2 = m_0 m'$ con $1 \le m' < m_0/2$.

Segue che:

$$m_0^2 m' p = (m_0 p)(m_0 m') = (x^2 + y^2)(x_1^2 + y_1^2) = (xx_1 + yy_1)^2 + (xy_1 - x_1 y)^2$$

Abbiamo anche:

$$xx_1 + yy_1 = x(x - cm_0) + y(y - dm_0) = (x^2 + y^2) - m_0(cx + dy) = m_0(p - cx - dy) = m_0 t$$

$$xy_1 - x_1 y = x(y - dm_0) - y(x - cm_0) = -m_0(xd - yc) = m_0 u$$

cioè esistono interi $t$ e $u$ tali che $m'p = t^2 + u^2$, con $1 \le m' < m_0/2 < m_0$. Questa è una contraddizione e conclude la dimostrazione. □

*Osservazione 9.2*

Nella dimostrazione del risultato precendente è fondamentale oltre alla Proposizione 3.2, anche la seguente identità già nota ad Eulero:

$$(a^2 + b^2)(c^2 + d^2) = (ac + bd)^2 + (ad - bc)^2 = (ac - bd)^2 + (ad + bc)^2 \qquad (9.2)$$

Al di là della verifica che si fa immediatamente svolgendo i calcoli, è interessante notare se si esce dall'insieme dei numeri interi per entrare in quello dei numeri complessi, la (9.2) è conseguenza del fatto se un numero complesso è prodotto di altri due allora il suo modulo è il prodotto dei relativi moduli, infatti:

$$(a^2 + b^2)(c^2 + d^2) = |(a + ib)(a - ib)(c + id)(c - id)|^2 =$$

$$|(a + ib)(c + id)|^2 \cdot |(a + ib)(c + id)|^2 = |(ac + bd) + i(ad - bc)|^2 = (ac + bd)^2 + (ad - bc)^2$$

$$|(a + ib)(c - id)|^2 \cdot |(a - ib)(c + id)|^2 = |(ac - bd) + i(ad + bc)|^2 = (ac - bd)^2 + (ad + bc)^2.$$

In realtà è proprio l'identità (9.2) che ha suggerito le regole di moltiplicazione alla base della costruzione dei numeri complessi.

*Definizione 9.3*-***Rappresentazione Primitiva***

Sia $n \in \mathbb{Z}$ con $n > 0$, si dice che $n$ ammette una rappresentazione primitiva come somma di due quadrati se esistono due interi $a$ e $b$ tali che $n = a^2 + b^2$ e $MCD(a,b) = 1$.

*Proposizione 9.4*

Se $p$ è un numero primo scrivibile come somma di due quadrati, allora la rappresentazione è primitiva ed è essenzialmente unica a parte il segno è l'ordine degli addendi.

Dimostrazione

Siano $a$ e $b$ interi tali che $p = a^2 + b^2$, allora posto $d = MCD(a,b)$ abbiamo che $d^2 \mid p$ da cui $d = 1$.

Se $p = 2 = a^2 + b^2$ allora $|a| < 2$ e $|b| < 2$ altrimenti $a^2 + b^2 > 2$ e $|ab| > 0$ altrimenti $a^2 + b^2 < p$

Quindi $|a| = 1$ e $|b| = 1$ che corrisponde alle quattro coppie: $(1,1),(1,-1),(-1,1),(-1,-1)$.

Se $p > 2$ allora per la Proposizione 9.1 ammette una rappresentazione come somma di due quadrati se e solo se $p \equiv 1 \pmod 4$. Se $p = a^2 + b^2 = c^2 + d^2$ dato che $p \neq 2$ possiamo escludere la possibilità che $|a| = |b|$ oppure che $|c| = |d|$, quindi senza perdere di generalità possiamo assumere $a$, $b$, $c$, $d$ interi positivi.

Allora risulta che:

$$p^2 = (a^2 + b^2)(c^2 + d^2) = (ac - bd)^2 + (ad + bc)^2 = (ac + bd)^2 + (ad - bc)^2$$

$$p(d^2 - b^2) = (a^2 + b^2)d^2 - (c^2 + d^2)b^2 = (ad - cb)(ad + cb) \equiv 0 \pmod p$$ quindi:

$(ad - cb) \equiv 0 \pmod p$ oppure $(ad + cb) \equiv 0 \pmod p$

essendo $0 < a,b,c,d < \sqrt{p}$ si ha che $ad - cb = 0$ oppure $ad + cb = p$.

Se $ad + cb = p$ allora essendo $p^2 = (a^2 + b^2)(c^2 + d^2) = (ac - bd)^2 + (ad + bc)^2 = (ac - bd)^2 + p^2$ segue che $ac = bd$ e dato che $MCD(a,b) = 1$ si ha che $a \mid d$ e quindi $d = ah$ da cui $c = hb$ ma essendo $1 = MCD(c,d) = h$ segue che $a = d$ e $c = b$.

Se $ad - cb = 0$ allora $a = c$ e $d = b$. □

*Corollario 9.5*

L'equazione diofantea $X^2 + Y^2 = p$ con $p$ primo dispari:

- **i.** non ammette soluzioni per $p \equiv 3 \pmod 4$
- **ii.** ammette otto soluzioni distinte nel caso $p \equiv 1 \pmod 4$, sia $p = a^2 + b^2$ una qualsiasi rappresentazione di $p$ come somma di due quadrati le otto soluzioni sono le seguenti: $(|a|,|b|)$, $(|a|,-|b|)$, $(-|a|,-|b|)$, $-(|a|,|b|)$, $(|b|,|a|)$, $(|b|,-|a|)$, $(-|b|,-|a|)$, $-(|b|,|a|)$.

Dimostrazione

Immediata conseguenza delle Proposizioni 9.4 e 9.1. □

Il caso più generale è quello relativo alla esistenza di una rappresentazione primitiva di un intero positivo come somma di due quadrati. Prima di analizzarlo in dettaglio conviene riassumere i risultati ottenuti nel caso di interi primi nel seguente teorema, questo ci consentirà di cogliere le analogie con il caso generale per cui enunceremo e dimostreremo un teorema perfettamente analogo.

*Teorema 9.6*

Sia $p$ un numero primo sono condizioni equivalenti le seguenti:

**i.** $p$ ammette una rappresentazione primitiva come somma di due quadrati.

**ii.** $p$ o $p/2$ è un intero dispari congruente ad 1 modulo 4 .

**iii.** La congruenza $X^2 \equiv -1 \pmod p$ è risolubile e ogni famiglia di soluzioni congruenti modulo $p$ determina una sola coppia di interi $(x,y)$ tali che:

  **1.** $kx \equiv y \pmod p$ con $k$ rappresentante qualsiasi della famiglia di soluzioni.

  **2.** $p = x^2 + y^2$ , $x > 0$ , $y > 0$, $MCD(x,y) = 1$

Dimostrazione

**i.** $\Leftrightarrow$ **ii.**

Vedi dimostrazione della Proposizione 9.1

**i.** $\Rightarrow$ **iii.**

Il caso $p = 2$ è immediato.

Sia quindi $p \equiv 1 \pmod 4$ dalla Proposizione 3.2(c) segue la risolubilità della congruenza $X^2 \equiv -1 \pmod p$, la quale come sappiamo ammette esattamente due soluzioni non congruenti modulo $p$. Sia ora $p = x^2 + y^2$ una rappresentazione primitiva con $x > 0$ e $y > 0$ allora $MCD(x,p) = MCD(y,p) = 1$ quindi siano, rispettivamente, $\overline{x}$ e $\overline{y}$ due inversi moltiplicativi modulo $p$. Abbiamo le due soluzioni non congruenti $k_1 = \overline{x}y$ , $k_2 = x\overline{y}$ di $X^2 \equiv -1 \pmod p$, inoltre $k_1 x \equiv y \pmod p$ e $k_2 y \equiv x \pmod p$.

**iii.** $\Rightarrow$ **i.**

Di immediata verifica dato che nella formulazione di iii. è contenuto anche i.. □

Prima di generalizzare il teorema precendente dimostriamo il seguente Lemma.

*Lemma 9.7*

Consideriamo gli interi $k,n \geq 1$ con $MCD(k,n) = 1$, allora la congruenza lineare in due indeterminate:

$kX \equiv Y \pmod n$ (9.3)

ammette sempre una soluzione $(x_0, y_0)$ con $0 < |x_0| < \sqrt{n}$ e $0 < |y_0| < \sqrt{n}$ .

Dimostrazione

Poniamo $N = \lfloor \sqrt{n} \rfloor + 1$ e consideriamo l'insieme:

$S = \{kx - y : 0 \leq x \leq N-1, 0 \leq y \leq N-1\}$ abbiamo che $|S| \leq N^2$, infatti posto $f(x,y) = kx - y$ , abbiamo che $S = f\left(\{1,2,\ldots,N-1\}^2\right)$ e $\left|\{1,2,\ldots,N-1\}^2\right| = N^2$.

Se $|S| < N^2$, allora devono esistere $(x_1,y_1) \neq (x_2,y_2)$ tali che $kx_1 - y_1 = kx_2 - y_2$ poniamo quindi $x_0 = x_2 - x_1$ e $y_0 = y_2 - y_1$ da cui $kx_0 = y_0$ e quindi in particolare anche $kx_0 = y_0 \pmod n$.

Se $|S| = N^2$, allora dato che $N^2 > n$ ed esistono solo $n$ classi di congruenza modulo $n$, due numeri $kx_1 - y_1$, $kx_2 - y_2$ generati da $(x_1,y_1) \neq (x_2,y_2)$ debbono appartenere alla stessa classe di congruenza modulo $n$, cioè $kx_1 - y_1 \equiv kx_2 - y_2 \pmod n$ quindi definendo $x_0$ e $y_0$ come nel caso precedente si ha il risultato. □

Siamo ora in grado di dimostrare la seguente generalizzazione del Teorema 9.6

*Teorema 9.8*

Sia $n$ un numero naturale sono condizioni equivalenti le seguenti:

**i.** $n$ ammette una rappresentazione primitiva come somma di due quadrati.

**ii.** $n$ o $n/2$ è un intero dispari ed è divisibile solo per numeri primi della forma $p \equiv 1 \pmod 4$.

**iii.** La congruenza $X^2 \equiv -1 \pmod n$ è risolubile e ogni famiglia di soluzioni congruenti modulo $n$ determina una ed una sola coppia di interi $(x, y)$ tale che:

  **1.** $kx \equiv y \pmod n$ con $k$ rappresentante qualsiasi della famiglia di soluzioni.

  **2.** $n = x^2 + y^2$, $x > 0$, $y > 0$, $MCD(x, y) = 1$

Dimostrazione

**i.** $\Rightarrow$ **ii.**

Se $n$ ammette una rappresentazione primitiva come somma di due quadrati allora esiste una coppia di numeri naturali $(x, y)$ tale che $n = x^2 + y^2$ e $MCD(x, y) = 1$. Sia quindi $p$ un numero primo tale che $p \mid n$, allora $p \nmid x \wedge p \nmid y$ altrimenti $p \mid MCD(x, y)$. Quindi esiste $k$ tale che $kx \equiv y \pmod n$. Allora $x^2 + y^2 \equiv x^2(1 + k^2) \pmod p$ quindi $p \mid (1 + k^2)$ cioè $p = 2$ oppure $p \equiv 1 \pmod 4$.

Se $n$ è pari allora si può avere $p = 2$, $x$ e $y$ sono dispari e quindi $n \equiv 2 \pmod 8$ da cui $n/2 \equiv 1 \pmod 4$ quindi $n/2$ è dispari. Se $p \neq 2$, allora $p \equiv 1 \pmod 4$ e $p \mid (n/2)$.

Se $n$ è dispari allora può solo essere $p \equiv 1 \pmod 4$.

**ii.** $\Rightarrow$ **iii.**

Applicazione del Teorema 7.7 segue che $X^2 \equiv -1 \pmod n$ è risolubile. Sia ora $k^2 \equiv -1 \pmod n$ con $1 \leq k \leq n$ si considerino $(x_0, y_0)$ come nel Lemma 9.7 poniamo $x = |x_0|$ e $y = |y_0|$.

- $\boxed{n = x^2 + y^2}$

$y_0 \equiv kx_0 \pmod n$ implica che ${y_0}^2 \equiv k^2 {x_0}^2 \pmod n$ quindi

$k^2(x^2 + y^2) = k^2({x_0}^2 + {y_0}^2) \equiv k^2 {x_0}^2 - {y_0}^2 \equiv 0 \pmod n$

siccome $MCD(k, n) = 1$ segue che $x^2 + y^2 \equiv 0 \pmod n$ ed essendo $0 < x^2 + y^2 < 2n$ segue che $x^2 + y^2 = n$.

- $\boxed{x > 0, y > 0}$

Se $x_0, y_0$ hanno segni concordi allora $kx \equiv y \pmod n$. Se invece $x_0, y_0$ hanno segni discordi ( ad esempio $x_0 < 0$ e $y_0 > 0$) allora poiché $k^2 \equiv -1 \pmod n$, abbiamo che $-ky_0 \equiv x_0 \pmod n$ ( infatti $kx_0 \equiv y_0 \pmod n$ da cui $x_0 \equiv -k^2 x_0 \equiv -ky_0 \pmod n$) quindi posto $x' = y$ e $y' = x$ abbiamo che $kx' \equiv y' \pmod n$.

- $\boxed{MCD(x, y) = 1}$

siano $f, g \in \mathbb{Z}$ tali che $k^2 = -1 + fn$ e $y = kx + gn$, allora:

$n = x^2 + y^2 = x^2 + (kx + gn)^2 = x^2(1 + k^2) + 2kxgn + g^2 n^2 = n[x(fx + kg) + g(kx + gn)] =$
$n[x(fx + kg) + gy]$ da cui $x(fx + kg) + gy = 1$ e per il teorema di Bezout $MCD(x, y) = 1$.

- $\boxed{\text{Biezione tra rappresentazioni primitive e soluzioni incongruenti di } X^2 \equiv -1 \pmod n}$

Se $(x_1, y_1)$ e $(x_2, y_2)$ sono due coppie di interi, per $k$ fissato, per cui valgono iii.1 e iii.2 allora abbiamo applicando la (9.2) che:

$n^2 = (x_1 x_2 + y_1 y_2)^2 + (x_1 y_2 - x_2 y_1)^2$ segue che $0 < x_1 x_2 + y_1 y_2 \leq n$

$x_1 x_2 + y_1 y_2 \equiv x_1 x_2 + (kx_1)(kx_2) \equiv (1 + k^2) x_1 x_2 \equiv 0 \pmod n$ da cui $n \mid (x_1 x_2 + y_1 y_2)$ e quindi $n = x_1 x_2 + y_1 y_2$ e $x_1 y_2 = x_2 y_1$ da cui essendo $MCD(x_1, y_1) = MCD(x_2, y_2) = 1$ segue che $x_1 = x_2$ e $y_1 = y_2$.

**iii.** $\Rightarrow$ **i.**

Di immediata verifica dato che nella formulazione di iii. è contenuto anche i.. □

*Corollario 9.9*

Sia $n$ un numero naturale tale $n = x^2 + y^2$ con $x, y$ naturali, allora se $p$ è un numero primo tale che $p \mid n$ e $p \equiv 3 \pmod 4$ allora $p \mid MCD(x, y)$.

Dimostrazione

Immediata conseguenza del Teorema 9.8, infatti altrimenti si contraddirebbe che i. implica ii. per il numero intero dotato di rappresentazione primitiva $(n / MCD(x, y))$. □

*Teorema 9.10*

Un numero naturale $n$ è la somma di due quadrati se e solo se ogni fattore primo $p$ di $n$, tale che $p \equiv 3 \pmod 4$ ha una potenza pari nella decomposizione di $n$ in fattori primi.

Dimostrazione

Sia $n = p_1^{e_1} \dots p_r^{e_r}$ e assumiamo che $e_j$ è pari quando $p \equiv 3 \pmod 4$. Allora $n = n_0 n_1^2$ dove $n_0 \geq 1$, $n_1 \geq 1$ e $n_0$ è il prodotto di primi distinti di cui al più uno è uguale a $2$ e gli altri sono congruenti ad $1$ modulo $4$. Applicando ripetutamente la formula (9.2) ad $n_0$ otteniamo la scrittura di $n_0$ e quindi $n$ come somma di due quadrati. Viceversa sia $n = x^2 + y^2$, la verifica è banale se $x = 0$ oppure $y = 0$, supponiamo $xy \neq 0$, e sia $d = MCD(x, y)$ allora $d^2 \mid n$. Poniamo $n = d^2 n'$, $x = dx'$, $y = dy'$, quindi $MCD(x', y') = 1$ e $n' = x'^2 + y'^2$. Applicando il Teorema 9.8 in particolare il fatto che i. implica iii. si ha il risultato. □

Al termine del prossimo paragrafo analizzeremo il numero di rappresentazioni distinte di un intero come somma di due quadrati.

## 10. *Anello degli interi di Gauss* $\mathbb{Z}(i)$

L'insieme di numeri complessi $\mathbb{Z}(i) = \{a + ib \mid a, b \in \mathbb{Z}\}$ considerato congiuntamente con le ordinarie operazioni di somma e moltiplicazione su $\mathbb{C}$ ha la struttura algebrica di anello (rif. [MA] o [A]).
$\mathbb{Z}(i)$ possiede proprietà che assomigliano a quelle dell'insieme $\mathbb{Z}$, in particolare esso possiede sotto opportune definizioni una proprietà di fattorizzazione unica analoga a quella di $\mathbb{Z}$.
Quello che vedremo in questo paragrafo nel caso particolare di $\mathbb{Z}(i)$ è generalizzato dalla teoria degli anelli euclidei, la quale assumendo che l'anello sia dotato di opportune proprietà, tali appunto da renderlo un "Anello Euclideo", generalizza l'algoritmo di Euclide, la fattorizzazioni unica e molte altre proprietà di $\mathbb{Z}$. Noi ci limiteremo allo studio $\mathbb{Z}(i)$ cercando di vedere sotto una differente luce i risultati del paragrafo precedente, il lettore interessato può divertirsi a tentare di generalizzare i concetti che vedremo in questo paragrafo o più semplicemente leggere una delle trattazioni fatte nei testi [AM], [A].
Di seguito per brevità chiameremo gli elementi di $\mathbb{Z}(i)$ interi, anziché interi Gaussiani.

*Definizione 10.1*
Un intero (Gaussiano) $\xi$ è divisibile dall'intero $\eta \neq 0$, oppure $\eta$ divide $\xi$, o infine $\eta$ è un divisore di $\xi$ abbreviato con $\eta \mid \xi$, se esiste un intero $\zeta$ tale che: $\xi = \eta\zeta$.

Derivano immediatamente dalla definizione le seguenti proprietà:
$\alpha \mid \beta \wedge \beta \mid \gamma \Rightarrow \alpha \mid \gamma$
$\alpha \mid \beta_1 \wedge \alpha \mid \beta_2 \Rightarrow \alpha \mid (\gamma_1\beta_1 + \gamma_2\beta_2)$

Un qualsiasi intero $\xi$ possiede i seguenti divisori banali: $1, \xi, -1, -\xi, i, i\xi, , -i, -i\xi$, questo suggerisce la generalizzazione del concetto di unità.

*Definizione 10.2*
Un intero $\varepsilon$ è detto una unità se $\varepsilon \mid \xi$ per ogni intero $\xi$.

*Proposizione 10.3*
Un intero $\varepsilon$ è una unità se e solo se è un divisore di $1$.

Dimostrazione
Se $\varepsilon$ è una unità in particolare $\varepsilon \mid 1$, viceversa se $\varepsilon \mid 1$ dato che $1 \mid \xi$, comunque sia $\xi$, segue che $\varepsilon \mid \xi$, comunque sia $\xi$, cioè $\varepsilon$ è una unità. □

*Definizione 10.4*
Gli interi $\xi$, $\eta$ sono associati se $\xi \mid \eta$ e $\eta \mid \xi$.

*Proposizione 10.5*
Due interi sono associati se e solo se differiscono l'uno dall'altro per un fattore moltiplicativo che è una unità di $\mathbb{Z}(i)$.

Dimostrazione
Siccome risulta $\xi = \alpha\eta$ e $\eta = \beta\xi$ si ha che $\alpha\beta = 1$ da cui $\alpha$ e $\beta$ sono unità. □

*Corollario 10.6*
Gli elementi $\mathbb{Z}(i)$ sono dotati di inverso moltiplicativo se e solo se sono unità.

Dimostrazione
Infatti se esiste $\alpha^{-1}$ allora $\alpha\alpha^{-1} = 1$ cioè $\alpha \mid 1$ e quindi $\alpha$ è una unità. Per l'implicazione inversa basta usare la Proposizione 10.3. □

*Definizione 10.7*
La norma di un intero $\xi = a + ib$ è definita come: $N(\xi) = a^2 + b^2$

Detto $\bar{\xi} = a - ib$ il complesso coniugato, allora $N(\xi) = \xi\bar{\xi} = |\xi|^2$, inoltre applicando la (9.2) si ha che posto $\eta = c + id$:
$N(\xi\eta) = N(\xi)N(\eta)$

*Proposizione 10.8*
Un intero $\varepsilon$ è una unità se e solo se ha norma $1$.

Dimostrazione
Se $\varepsilon$ è una unità, allora $\varepsilon \mid 1$, quindi esiste un intero $\eta$ tale che $1 = \varepsilon\eta$, quindi $1 = N(\varepsilon)N(\eta)$ da cui $N(\varepsilon) = 1$. Viceversa se $1 = N(\varepsilon) = N(a + ib) = a^2 + b^2 = (a + ib)(a - ib)$ segue che $\varepsilon \mid 1$ e applicando la Proposizione 10.3 segue il risultato. □

*Proposizione 10.9*
Le unita di $\mathbb{Z}(i)$ sono $1$, $-1$, $i$, $-i$.

Dimostrazione
Le uniche soluzioni di $a^2 + b^2 = 1$, sono $a = \pm 1, b = 0$ e $a = 0$, $b = \pm 1$. □

*Osservazione 10.10*
Gli associati dell'intero $\xi$ sono quindi: $\xi$, $i\xi$, $-\xi$, $-\xi$.

*Definizione 10.11-**Primi di $\mathbb{Z}$(i)***

Un primo è un intero $\pi$ di $\mathbb{Z}(i)$ che non una unità ed è diverso da 0, divisibile solo per i suoi associati e per le unità.

*Proposizione 10.12*
Sia $\xi$ un intero con $N(\xi) = p$ primo di $\mathbb{Z}$, allora $\xi$ è un primo di $\mathbb{Z}(i)$.

Dimostrazione
Infatti se $\eta \mid \xi$ e $N(\eta) \neq 1$, allora $\xi = \alpha\eta$ e $N(\xi) = N(\alpha\eta) = N(\alpha)N(\eta) = p$, quindi $N(\eta) \mid p$ cioè $N(\eta) = p$, da cui segue che $N(\alpha) = 1$ cioè $\eta$ è un associato di $\xi$. □

*Osservazione 10.13*
Ritroviamo sotto una differente veste i concetti del paragrafo precedente dove abbiamo visto che $p = a^2 + b^2 \Rightarrow (p \equiv 1 \pmod 4) \vee (p = 2)$, inoltre le diverse rappresentazioni come somme di due quadrati del numero primo $p$ di $\mathbb{Z}$ corrispondono ai diversi interi associati a $\pi = a + ib$ al numero primo di $\mathbb{Z}(i)$, cioè $\pi$, $i\pi$, $-\pi$, $-i\pi$.

*Proposizione 10.14*

Ogni primo $q$ di $\mathbb{Z}$ della forma $q \equiv 3 \pmod 4$ è anche un primo di $\mathbb{Z}(i)$.

Dimostrazione

Sia $\xi$ tale che $\xi \mid q$ e $N(\xi) \neq 1$, allora $q = \alpha\xi$ e quindi $N(q) = N(\alpha)N(\xi) = q^2$ da cui $N(\xi) \mid q^2$, dato che non può essere $N(\xi) = q$, dato che $q \equiv 3 \pmod 4$, si ha che $N(\xi) = q^2$ cioè $N(\alpha) = 1$ e $\xi$ associato di $q$.

*Osservazione 10.15*

Dalle precedenti due Proposizioni segue che non tutti i primi di $\mathbb{Z}$ sono primi di $\mathbb{Z}(i)$. Questo fatto è essenzialmente equivalente alla Proposizione 9.1.

*Teorema 10.16*

I primi di $\mathbb{Z}(i)$ sono:

- **i.** $q \in \mathbb{Z}$ primo tale che $q \equiv 3 \pmod 4$.
- **ii.** $\xi \in \mathbb{Z}(i)$ tali che $N(\xi) = p$ con $p$ primo di $\mathbb{Z}$ e $p \equiv 1 \pmod 4 \vee p = 2$.

Dimostrazione

Sia $\pi \in \mathbb{Z}(i)$ un primo senza perdere di generalità possiamo assumere considerando eventualmente un suo associato al posto di $\pi$, che $\pi = a + ib$ con $a > 0$. Consideriamo quindi $N(\pi) = N(a + ib) = a^2 + b^2 \in \mathbb{Z}$ allora definiamo $q = MCD(a,b)$, $a' = a/q$ e $b' = b/q$ allora risulta che $\pi = (q + i0)(a' + ib')$ e per il Teorema 9.8 si hanno le seguenti due possibilità:

- $b' = 0$ da cui segue che $a' = 1$ perché $MCD(a',b') = 1$, e $\pi = q$ di conseguenza $q \in \mathbb{Z}$ primo e $q \equiv 3 \pmod 4$ per non contraddire la primalità di $\pi$.
- $b' \neq 0$ allora essendo $\pi$ primo si ha $q = 1$, e si presentano due sottocasi (Teorema 9.8):
  - Se $2 \mid N(\pi)$ segue per la primalità di $\pi$ che $N(\pi) = 2$.
  - Se invece $2 \nmid N(\pi)$ allora $N(\pi)$ è il prodotto di primi di $\mathbb{Z}$ del tipo $1 \pmod 4$, e data la primalità di $\pi$ segue che $N(\pi) = p$ con $p \equiv 1 \pmod 4$.

□

Stabilito quali sono i primi di $\mathbb{Z}(i)$ procediamo alla dimostrazione di una serie di Teoremi propedeutici alla dimostrazione del Teorema fondamentale dell'aritmetica per gli interi Gaussiani.

*Teorema 10.17*

Ogni intero diverso da $0$ e che non sia una unità è divisibile per un numero primo.

Dimostrazione

Se $\gamma$ è un intero non primo, allora:

$\gamma = \alpha_1\beta_1$, $N(\alpha_1) > 1$, $N(\beta_1) > 1$, $N(\gamma) = N(\alpha_1)N(\beta_1)$, e quindi $1 < N(\alpha_1) < N(\gamma)$.

Se $\alpha_1$ non è un primo, allora:

$\alpha_1 = \alpha_2\beta_2$, $N(\alpha_2) > 1$, $N(\beta_2) > 1$, $N(\gamma) = N(\alpha_2)N(\beta_2)$, e quindi $1 < N(\alpha_2) < N(\alpha_1)$.

Questo processo può essere continuato fino a che $\alpha_r$ non è un primo. Siccome $N(\gamma)$, $N(\alpha_1)$, $N(\alpha_2)$, ... è una successione decrescente di interi positivi, si deve prima o poi raggiungere un primo $\alpha_r$ ,e ammesso che sia il primo numero primo ad apparire nella sequenza. Allora:

$\gamma = \beta_1\alpha_1 = \beta_1\beta_2\alpha_2 = \ldots = \beta_1\beta_2\beta_2\ldots\beta_r\alpha_r$ e quindi $\alpha_r \mid \gamma$. □

*Teorema 10.18*

Ogni intero diverso da $0$ e che non sia una unità è il prodotto di primi.

Dimostrazione
Se $\gamma$ diverso da $0$ e non è una unità allora per il Teorema 10.17 è divisibile per un primo $\pi_1$. Quindi $\gamma = \pi_1\gamma_1$, $N(\gamma_1) < N(\gamma)$. Se $\gamma_1$ non è una unità allora $\gamma_1 = \pi_2\gamma_2$, $N(\gamma_2) < N(\gamma_1)$. Continuando questo processo otteniamo una successione decrescente di interi positivi $N(\gamma)$, $N(\gamma_1)$, $N(\gamma_2)$,.... Quindi $N(\gamma_r) = 1$ per $r$ opportuno, e $\gamma_r$ è una unità posto quindi $\pi'_r = \pi_r\gamma_r$, si ha che $\pi_r$ e $\pi'_r$ sono associati e possiamo scrivere $\gamma = \pi_1\pi_2 ... \pi'_r$. □

*Teorema 10.19-**Algorimo della divisione in $\mathbb{Z}(i)$***

Dati due interi $\alpha$, $\beta$, dei quali $\beta \neq 0$ esiste un intero $\kappa$ tale che:
$\alpha = \kappa\beta + \rho_1$, $N(\rho_1) < N(\beta)$.

Dimostrazione
Siccome $\beta \neq 0$, abbiamo $\alpha / \beta = R + Si$ dove $R$ ed $S$ sono numeri razionali. Possiamo trovare sue interi tali che $|R - x| \leq 1/2$, $|S - y| \leq 1/2$. Ad esempio basta prendere $x = \lfloor R \rfloor$ se $R \leq \lfloor R \rfloor + 1/2$ e $x = \lfloor R + 1 \rfloor$ se $R > \lfloor R \rfloor + 1/2$ ed analogamente per $y$ scambiando $R$ con $S$. Abbiamo quindi che:
$$|\alpha / \beta - (x + iy)| = |(R - x) + i(S - y)| = \sqrt{(R - x)^2 + (S - y)^2} \leq 1/\sqrt{2}$$
Se quindi prendiamo $\kappa = x + iy$, $\rho_1 = \alpha - \kappa\beta$, allora abbiamo $|\rho_1| = |\alpha - \kappa\beta| \leq |\beta| / \sqrt{2}$ da cui quadrando $N(\rho_1) = N(\alpha - \kappa\beta) \leq N(\beta)/2$. □

*Definizione 10.20-**Massimo comun divisore in $\mathbb{Z}(i)$***

Se $\zeta$ è un divisore comune di $\alpha$ e $\beta$, ed ogni divisore comune di $\alpha$ e $\beta$ è un divisore di $\zeta$, allora chiamiamo $\zeta$ *un massimo comun divisore* di $\alpha$ e $\beta$.
A differenza di $\mathbb{Z}$, in $\mathbb{Z}(i)$ il massimo comun divisore non è unico in quanto ogni associato ad un massimo comun divisore è un massimo comun divisore, quindi quando scriviamo $\zeta = MCD(\alpha, \beta)$ intendiamo che $\zeta$ è un $MCD$ di $\alpha$ e $\beta$. Quando il $MCD$ di $\alpha$ e $\beta$ esiste è unico a meno di associati, infatti due $MCD$ distinti debbono dividersi l'un l'altro e quindi sono associati.

*Teorema 10.21-**Algoritmo di Euclide in $\mathbb{Z}(i)$***

Siamo $\alpha$ e $\beta$ due interi con $N(\alpha) > N(\beta)$ allora esiste $MCD(\alpha, \beta)$.

Dimostrazione
Chiaramente abbiamo $N(\alpha) > 0$ e se $\beta = 0$ allora $\zeta = \alpha$. Supponiamo quindi $N(\alpha) \geq N(\beta) > 0$. Dividendo $\alpha$ per $\beta$ otteniamo: $\alpha = \kappa_1\beta + \rho_1$ con $0 \leq N(\rho_1) < N(\beta)$. Se $\rho_1 \neq 0$ allora possiamo ripetere il processo ottenendo $\beta = \kappa_2\rho_1 + \rho_2$ dove $0 \leq N(\rho_2) < N(\rho_1)$. Se $\rho_2 \neq 0$, $\rho_1 = \kappa_3\rho_2 + \rho_3$, dove $0 \leq N(\rho_3) < N(\rho_2)$ e così via. La successione di interi non negativi $N(\beta)$, $N(\rho_1)$, $N(\rho_2)$,...è decrescente e quindi per $n$ opportuno avremo $\rho_{n+1} = 0$. Gli ultimi due passi del processo sono: $\rho_{n-2} = \kappa_n\rho_{n-1} + \rho_n$ dove $0 < \rho_n < \rho_{n-1}$, $\rho_{n-1} = \kappa_{n+1}\rho_n$. Verifichiamo che risulta quindi che $\rho_n = MCD(\alpha, \beta)$. Procedento a ritroso $\rho_n \mid \rho_{n-1} \Rightarrow \rho_n \mid \rho_{n-2} \Rightarrow ... \Rightarrow \rho_n \mid \alpha \wedge \rho_n \mid \beta$ cioè $\rho_n$ è un divisore comune. Se $\zeta \mid \alpha \wedge \zeta \mid \beta \Rightarrow \zeta \mid \rho_1 \Rightarrow \zeta \mid \rho_2 \Rightarrow ... \Rightarrow \zeta \mid \rho_n$. Segue per la definizione di massimo comun divisore che $\rho_n = MCD(\alpha, \beta)$. Se $\zeta$ è un massimo comun divisore allora risulta in particolare che $\zeta \mid \rho_n$ e quindi $\zeta$ è associato a $\rho_n$. □

*Teorema 10.22*

Se $MCD(\alpha, \beta) = 1$ e $\alpha \mid \beta\gamma$, allora $\alpha \mid \gamma$.

Dimostrazione
Pensando di applicare l'algoritmo di Euclide giungiamo a $\rho_n = MCD(\alpha, \beta)$ moltiplicando ambo i membri dell'insieme di equazioni che costituiscono l'algoritmo per $\gamma$, abbiamo che $\gamma\rho_n = MCD(\alpha\gamma, \beta\gamma)$ per ipotesi $\alpha \mid \beta\gamma$ e inoltre $\alpha \mid \alpha\gamma$ segue che $\alpha \mid \gamma\rho_n$. Dato che $\rho_n$ è associato ad $1$ si ha infine che $\alpha \mid \gamma$. □

*Teorema 10.23*
Se $\pi$ è un primo e $\pi \mid \alpha\beta$, allora $\pi \mid \alpha$ o $\pi \mid \beta$.

Dimostrazione
Sia $\mu = MCD(\pi, \alpha)$ allora se $\mu$ è associato ad $1$, abbiamo che per il Teorema 10.22 $\pi \mid \beta$, se invece $\mu$ è associato ad $\pi$ segue che $\pi \mid \alpha$. □

*Teorema 10.24-**Teorema fondamentale per gli interi Gaussiani.***
La scomposizione di un intero diverso da $0$ o da una unità come prodotto di primi è unica, a meno dell'ordine dei primi, la presenza di unità e ambiguità tra primi associati.

Dimostrazione
Assumiamo per assurdo che il teorema sia falso e sia $\gamma$ un intero per cui il teorema è falso con $N(\gamma)$ è minimo, nel senso che ogni intero $\alpha$ con $N(\alpha) < N(\gamma)$ il teorema è vero, chiaramente dato che $\gamma$ non è una unità abbiamo $N(\gamma) > 1$. Allora esistono almeno due fattorizzazioni $\gamma = \pi_1 ... \pi_r = \pi'_1 ... \pi'_2$, allora $\pi_1 \mid \pi'_1 ... \pi'_s$, quindi per il Teorema 10.23 deve dividere un $\pi'_j$ senza perdere di generalità possiamo assumere che sia $\pi'_1$. Segue che $\pi_1$ e $\pi'_1$ primi sono associati, $N(\gamma / \pi_1) < N(\gamma)$ e $\gamma / \pi_1$ ammette due fattorizzazioni da cui l'assurdo. □

Siamo ora in grado di affrontare il calcolo del numero di rappresentazioni di un intero come somma di due quadrati.

*Definizione 10.25-**La funzione $r(n)$***

Definiamo $r(n)$ come il numero di rappresentazioni di $n$ nella forma $n = A^2 + B^2$ dove $A$ e $B$ sono interi. Consideriamo distinte due rappresentazioni anche se loro differiscono banalmente, cioè relativamente al segno e all'ordine di $A$ e $B$.

Quindi:
$0 = 0^2 + 0^2$, $r(0) = 1$.
$1 = (\pm 1)^2 + 0^2 = 0^2 + (\pm 1)^2$, $r(1) = 4$.
$r(p) = 8$, se $p$ è un primo tale che $p \equiv 1 \pmod 4$.
$r(q) = 0$, se $q$ è un primo tale che $q \equiv 3 \pmod 4$.
Dimostriamo ora il seguente teorema generale.

*Teorema 10.26*
$r(n) = 4(d_1(n) - d_3(n))$
dove:
$d_1(n)$ è il numero di divisori di $n$ congruenti ad $1$ modulo $4$.
$d_3(n)$ è il numero di divisori di $n$ congruenti a $3$ modulo 4.

Dimostrazione
Possiamo scrivere $n$ come segue:

$$n = 2^g \prod_{r=1}^{R} p^{e_r} \prod_{s=1}^{S} q^{f_s} = [(1+i)(1-i)]^g \prod_{r=1}^{R} [(a_r + ib_r)(a_r - ib_r)]^{e_r} \prod_{s=1}^{S} q^{f_s} \qquad (10.1)$$

dove:

- Per $r = 1,\dots,R$ $p_r \equiv 1 \pmod 4$ è un primo di $\mathbb{Z}$, $a_r$ e $b_r$ sono interi positivi distinti tali che $a_r > b_r$, $p_r = a_r^2 + b_r^2$, e questa rappresentazione di $p_r$ è unica. Inoltre $a_r + ib_r$, $a_r - ib_r$ sono primi di $\mathbb{Z}(i)$.
- Per $s = 1,\dots,S$, $q_s \equiv 3 \pmod 4$ è un primo di $\mathbb{Z}$ ed un primo di $\mathbb{Z}(i)$.
- $1+i$, $1-i$ sono primi di $\mathbb{Z}(i)$.

Quindi se pensiamo $n \in \mathbb{Z}(i)$ per il Teorema 10.24 la scrittura in forma complessa della (10.1) è unica a meno di unità e primi tra loro associati.
Se ora $n = A^2 + B^2$ è una generica rappresentazione come somma di due quadrati abbiamo che $n = (A+iB)(A-iB)$ considerando di nuovo $n \in \mathbb{Z}(i)$ deve risultare che $A+iB$ possiede una frazione dei fattori primi complessi di $n$, quindi possiamo assumere la forma generale.

$$A + iB = \varepsilon(1+i)^{g'}(1-i)^{g''}\prod_{r=1}^{R}(a_r + ib_r)^{e'_r}(a_r - ib_r)^{e''_r}\prod_{s=1}^{S}q^{f'_s} \qquad (10.2)$$

con $\varepsilon$ un unità di $\mathbb{Z}(i)$ costituita come il prodotto di tutte le potenze di unità, che sono a loro volta unità e renderebbero ambigua la rappresentazione come somma di due quadrati di $2$ e $p_r$ con $r = 1,\dots,R$. Prendendo il complesso coniugato abbiamo quindi che:

$$A - iB = \bar{\varepsilon}(1-i)^{g'}(1+i)^{g''}\prod_{r=1}^{R}(a_r - ib_r)^{e'_r}(a_r + ib_r)^{e''_r}\prod_{s=1}^{S}q^{f'_s}$$

imponendo l'uguaglianza del prodotto ad $n$ deve risultare $g' + g'' = g$, $e'_r + e''_r = e_r$, $2f'_s = f_s$.
Ritroviamo di nuovo che $2 \mid f_s$ per $s = 1,\dots,S$, quindi assumiamo queste condizioni per cui $r(n) \neq 0$. Essendo $\frac{1-i}{1+i} = -i$ e $\frac{1+i}{1-i} = i$ segue che:

$$A + iB = \varepsilon'(1-i)^{g}\prod_{r=1}^{R}(a_r + ib_r)^{e'_r}(a_r - ib_r)^{e''_r}\prod_{s=1}^{S}q^{f_s/2} \qquad (10.3)$$

$$A - iB = \bar{\varepsilon}'(1+i)^{g}\prod_{r=1}^{R}(a_r - ib_r)^{e'_r}(a_r + ib_r)^{e''_r}\prod_{s=1}^{S}q^{f_s/2} \qquad (10.4)$$

con $\varepsilon' = \varepsilon i^{g'}$ unità.
Pertanto le possibili rappresentazioni di $n$ come somma di due quadrati corrispondono alle possibili scelte indipendenti della unità $\varepsilon'$ e delle coppie $(e'_r, e''_r)$ tali che $e'_r + e''_r = e_r$ con $r = 1,\dots,R$. Quindi abbiamo che:

$$r(n) = 4 \cdot \prod_{r=1}^{R}(e_r + 1) \qquad (10.5)$$

Per contare i divisori dispari di $n$ basta notare che questi sono gli addendi del prodotto:

$$\prod_{r=1}^{R}\left(1 + p_r + \dots + p_r^{e_r}\right)\prod_{s=1}^{S}\left(1 + q_s + \dots + q_s^{f_s}\right)$$

un divisore è della forma $4m+1$ se contiene un numero pari di fattori $q$, mentre un divisore è della forma $4m-1$ se contiene un numero dispari di fattori $q$. Sostituendo 1 per $p_r$ e $-1$ per $q_s$ nel prodotto abbiamo che si ottiene la quantità $d_1(n) - d_3(n)$, quindi:

$$d_1(n) - d_3(n) = \prod_{r=1}^{R}\left(1 + 1 + \dots + 1^{e_r}\right)\prod_{s=1}^{S}\left(1 + (-1) + \dots + (-1)^{f_s}\right) = \prod_{r=1}^{R}(1 + e_r)\prod_{s=1}^{S}\frac{1 + (-1)^{f_s}}{2}$$

dato che il secondo fattore è nullo proprio quando uno qualsiasi degli $f_s$ non è pari segue dalla (10.5) la formula generale:

$$r(n) = 4\left(d_1(n) - d_3(n)\right)$$ □

*Corollario 10.27*

Se $n = 2^g \prod_{r=1}^{R} p^{e_r} \prod_{s=1}^{S} q^{f_s}$, allora $r(n) = 4\prod_{r=1}^{R}(1 + e_r)\prod_{s=1}^{S}\frac{1 + (-1)^{f_s}}{2}$.

Dimostrazione
Vedi Teorema 10.26. □

*Corollario 10.28*

La generica rappresentazione primitiva di $n = 2^g \prod_{r=1}^{R} p^{e_r} = A^2 + B^2$ è data da:

$A + iB = \varepsilon'(1-i)^g \prod_{r=1}^{R} (a_r + ib_r)^{\lambda'_r e_r} (a_r - ib_r)^{\lambda''_r e_r}$ con $(\lambda'_r, \lambda''_r) \in \{0,1\}^2$ e $\lambda'_r + \lambda''_r = 1$ per $1 \le r \le R$ e $g = 0$ se $n$ è dispari e $g = 1$ se $n$ è pari. Il numero di rappresentazioni primitive è $2^R$.

Dimostrazione

Per il Teorema 9.8 $g \in \{0,1\}$ a seconda che $n$ sia pari o dispari. Si consideri la formula (10.3) e riscriviamo il generico fattore come segue:

$$(a_r + ib_r)^{e'_r}(a_r - ib_r)^{e''_r} = \begin{cases} (a_r^2 + b_r^2)^{e'_r}(a_r - ib_r)^{e''_r - e'_r} = p_r^{e'_r}(a_r - ib_r)^{e''_r - e'_r} & e'_r \le e''_r \\ (a_r + ib_r)^{e'_r - e''_r}(a_r^2 + b_r^2)^{e''_r} = (a_r + ib_r)^{e'_r - e''_r} p_r^{e''_r} & e'_r > e''_r \end{cases}$$

Segue che $p^{\min(e'_r, e''_r)} \mid MCD(A,B)$ quindi $MCD(A,B) = 1$ implica che $e'_r = 0 \wedge e''_r = e_r$ oppure $e'_r = e_r \wedge e''_r = e_r$, da cui $e'_r = \lambda'_r e_r$ e $e''_r = \lambda''_r e_r$ con $\lambda'_r + \lambda''_r = 1$. Segue immediatamente che le possibili scelte delle $R$ coppie $(\lambda'_r, \lambda''_r) \in \{0,1\}^2$ nell'insieme con il vincolo $\lambda'_r + \lambda''_r = 1$ è $2^R$. □

## 11. L'equazione Diofantea $X^2 + Y^2 = cZ^2$ in $\mathbb{Z}$

Per $c > 0$ intero assegnato ci occupiamo di determinare tutte le terne di interi $(x, y, z)$ soluzioni della equazione diofantea:

$$X^2 + Y^2 = cZ^2 \qquad (11.1)$$

Possiamo limitarci alla ricerca di soluzioni primitive per cui:

$$MCD(x, y, z) = 1 \qquad (11.2)$$

Siamo interessati alle soluzioni non banali cioè diverse dalla terna $(0,0,0)$. Quindi in generale si dira la (1.1) risolvibile se ammette soluzioni differenti da quella banale.
Conviene iniziare il nostro studio ricordando il caso particolare della terne pitagoriche, corrispondenti al caso $c = 1$. La letteratura on-line e non sull'argomento è veramente abbondante, infatti un semplice "query" su un qualsiasi motore di ricerca da luogo a migliaia di risultati. Quindi limitiamoci ad enunciare il teorema generale per le formule delle terne pitagoriche primitive, la dimostrazione sarà un caso particolare del risultato più generale che otterremo per la (11.1).

*Teorema 11.1*

Se $m$ e $n$ sono interi tali che $m > n > 0$, $MCD(m,n) = 1$, $m$, $n$ di differente parità, allora le terne $(s, t, r)$, date da:

$$\begin{cases} s = 2mn \\ t = m^2 - n^2 \\ r = m^2 + n^2 \end{cases} \qquad (11.3)$$

è una soluzione primitiva della (11.2) per $c = 1$. Questo stabilisce una corrispondenza biunivoca tra l'insieme di coppie $(m, n)$ che soddisfano le condizioni sopra definite e l'insieme delle soluzioni primitive della (11.1) con $c = 1$ e $x$ pari.

Dimostrazione
Risulta un caso particolare della dimostrazione che daremo per i Teorema 11.6, 12.1 e 13.1
. □

Quindi la generica terna soluzione della (11.1) è una generalizzazione delle terne pitagoriche, che in analogia potremmo chiamare terna "quasi" pitagorica. Tuttavia un termine analogo potrebbe essere applicato alle soluzioni delle equazione più generale $aX^2 + bY^2 = cZ^2$ con $a$, $b$, $c$ interi positivi di cui la (11.1) è a sua volta un caso particolare. Nella discussione della Proposizione 7.8 sono state definite le condizioni necessarie e sufficienti per l'esistenza di una soluzione che corrisponde ad un triangolo di lati multipli interi dei numeri, in generale irrazionali, corrispondenti alla radici quadrate dei coefficienti dell'equazione: $x\sqrt{a}$, $y\sqrt{b}$, $z\sqrt{c}$.
Relativamente all'uso del nome di terne quasi pitagoriche si trova nella letteratura reperibile on-line che vengono così chiamate anche le terne corrispondenti ai triangoli di lati interi $x$, $y$, $z$ con l'angolo opposto al lato $z$ di $120°$. Per il teorema di Carnot risulta soddisfatta l'equazione diofantea $Z^2 = X^2 + XY + Y^2$.
Quindi senza formalizzarsi troppo sul nome da dare alle soluzioni della (11.1), anticipiamo che uno dei risultati principale di questo paragrafo consisterà proprio nel verificare che ogni soluzione può essere generata da una terna pitagorica tramite l'applicazione di una opportuna trasformazione di coordinate, e che questo risultato risulta immediato da verificare in $\mathbb{Z}(i)$ utilizzando la fattorizzazione unica.

*Teorema 11.2*

La (11.1) è risolubile se e solo se $c$ è scrivibile come somma di due quadrati.

Dimostrazione
La condizione è sufficiente in quanto se esistono due interi $u$ e $v$ tali che $c = u^2 + v^2$ allora la terna $(u, v, 1)$ è soluzione della 11.1. Verifichiamo che tale condizione è necessaria, se $c$ non è rappresentabile come somma di due quadrati, allora per il Teorema 9.10 $c$ ammette fattori primi del tipo $3 \pmod 4$ con esponente dispari, quindi se per la (11.1) ammettesse una soluzione $(x,y,z)$ non banale, ma dato che tutti gli esponenti dei fattori primi di $z^2$ sono pari ne risulterebbe che $cz^2$ è rappresentabile come somma di due quadrati e possiete fattori primi del tipo 3 (mod 4) con esponente dispari da cui l'assurdo. □

*Osservazione 11.3*

Un sotto insieme particolare della soluzioni primitive è definito dalla condizione più restrittiva:

$$d_3 = MCD(x,y) = 1 \qquad (11.4)$$

Per il Teorema 9.8 se la (11.1) ammette soluzioni che soddisfano la (11.4) allora $cz^2$, e quindi $c$, ha come divisori primi 2 o primi del tipo $1 \pmod 4$.

In generale per un intero positivo $c$ assegnato una generica soluzione primitiva $(x,y,z)$ della (11.1) può essere ricondotta ad risoluzione della (11.1) che rispetta la (11.4) definita dal parametro:

$$c' = c / {d_3}^2 \qquad (11.5)$$

sia $(x',y',z')$ con $MCD(x',y') = 1$ è una soluzione della

$$X^2 + Y^2 = c'Z^2 \qquad (11.6)$$

allora posto $(d_3x', d_3y', z')/MCD(d_3, z')$ è una soluzione della (11.1). Viceversa se $(x,y,z)$ è una soluzione primitiva della (11.1) allora $(x/d_3, y/d_3, z)$ allora $MCD(x/d_3, y/d_3) = 1$ inoltre $MCD(d_3, z) = MCD(x,y,z) = 1$. Quindi abbiamo dimostrato la seguente proposizione:

*Proposizione 11.4*

L'insieme delle soluzione primitive della (11.1) è costituito dall'unione dell'insieme delle soluzioni che soddisfano la (11.4) e degli insiemi analoghi ottenuti per ciascun divisore quadratico $d_3$ di $c$, sostituendo nella (11.1) il parametro $c$ con il parametro definito dalla (11.5) e con la clausola aggiuntiva $MCD(d_3, z) = 1$.

Dimostrazione

Vedi Osservazione 11.3. □

*Corollario 11.5*

Se $c$ non possiede divisori quadratici, la (11.1) se risolubile ammette solo soluzioni che soddisfano la (11.4).

Dimostrazione

Immediata conseguenza della Proposizione 11.4 o direttamente della (11.5). □

Essendoci ridotti sempre al caso di soluzioni della (11.1) per cui vale la (11.4), le formule risolutive generali della (11.1) sono definite dal seguente teorema.

*Teorema 11.6*

La generica soluzione primitiva $(x,y,z)$ della (11.1) è definita dalle seguenti equazioni:

${d_3}^2 \mid c$ e $u > v > 0$ con $MCD(u,v) = 1$ tali che $c/{d_3}^2 = 2^g(u^2 + v^2)$ e $g \in \{0,1\}$

$r^2 = s^2 + t^2$ terna pitagorica primitiva definita dalle formule (11.3) con la condizione aggiuntiva $MCD(r, d_3) = 1$.

$(x,y,z) \in \{(j_1\hat{x}, j_2\hat{y}, j_3\hat{z}), (j_1\hat{y}, j_2\hat{x}, j_3\hat{z}) \mid j_k = -1,1; k = 1,2,3\}$

Se $g = 0$ $$\begin{cases} \hat{x}/d_3 = (tu - sv) \\ \hat{y}/d_3 = (su + tv) \end{cases} \qquad (11.7)$$

Se $g = 1$ $$\begin{cases} \hat{x}/d_3 = (s+t)u - (s-t)v \\ \hat{y}/d_3 = (s-t)u + (s+t)v \end{cases} \qquad (11.8)$$

$$\hat{z} = r \qquad (11.9)$$

Dimostrazione
Svolgiamo nel dettaglio la dimostrazione per $d_3 = 1$ di cui tutti gli altri casi saranno immediata conseguenza ripedendo tutti i ragionamenti per $c/{d_3}^2$. Sia quindi $d_3 = 1$, se esistono $(x,y,z)$ tali che $MCD(x,y) = 1$ e $cz^2 = x^2 + y^2$, cioè $cz^2$ ammette una rappresentazione primitiva, allora per il Corollario 10.28 $z$ deve essere dispari altrimenti 2 avrebbe un esponente $g \geq 2$. Fattorizzando in $\mathbb{Z}(i)$ il numero $cz^2$ abbiamo che:

$cz^2 = (u+iv)(u-iv)(m+in)^2(m-in)^2$

dove:

$|z| = r = m^2 + n^2$ con $MCD(m,n) = 1$ e $m > n$.

$c = 2^g(u^2 + v^2)$ con $MCD(u,v) = 1$ e $u > v$.

e applicando il Corollario 10.28 risulta quindi:

$x + iy = i^r(1-i)^g(u+iv)(m+in)^2 = i^r(\hat{x} + i\hat{y})$ con $r = 0,1,2,3$ e $g = 0$ se $c$ è dispari, $g = 1$ se $c$ è pari.

Posto quindi:

$s = m^2 - n^2$

$t = 2mn$

svolgendo i calcoli si trovano le formule dell'asserto del teorema. □

*Osservazione 11.7*

Analizziamo i risultati del precendente teorema in analogia al caso delle terne pitagoriche per cui si ha $c = 1$. Risulta evidente che ogni rappresentazione di $c$ come somma di due quadrati, sia essa primitiva o non, genera un ramo (o famiglia o insieme) di soluzioni della (11.1). Possiamo trascurare, considerandole equivalenti, le rappresentazione di $c$ che differiscono solo per segno o per ordine delle componenti, dato che danno luogo a soluzioni che differiscono tra loro, a parità di terna pitagorica, solo per segno e ordine delle componenti. Pertanto, il caso $c = 1$ da luogo ad un solo ramo di soluzioni, corrispondenti alla unica rappresentazione $c = 1^2 + 0^2$ e lo stesso vale per il caso $c = 2 = 1^2 + 1^2$. Quindi nel senso appena precisato la (11.1) presenta tanti rami di soluzioni quante sono le rappresentazioni di $c$ come somma di due quadrati, il cui numero è un ottavo del numero definito formula (10.5) quando $c > 2$.

Una interpretazione geometrica degli argomenti svolti in questo paragrafo, si ottiene tenendo presente il fatto che indirettamente abbiamo risolto i seguente problema:

> Assegnata una terna pitagorica primitiva $(s,t,r)$ ed un numero intero positivo $c$ che non è un quadrato perfetto, determinare le rotazioni $\varphi$ del raggio vettore $r\sqrt{c}$ affinché le sue proiezioni sugli assi siano segmenti di lunghezza intera.

## 12. L'equazione Diofantea $X^2 + Y^2 = Z^l$ in $\mathbb{Z}$

Per $l > 1$ intero assegnato ci occupiamo di determinare tutte le terne di interi $(x, y, z)$ soluzioni primitive della equazione diofantea:

$$X^2 + Y^2 = Z^l \qquad (12.1)$$

In questo caso particolare se $x^2 + y^2 = z^l$:

$$MCD(x, y, z) = 1 \Leftrightarrow MCD(x, y) = 1 \Leftrightarrow MCD(x, z) = 1 \Leftrightarrow MCD(y, z) = 1 \qquad (12.2)$$

Quindi se $(x, y, z)$ è una soluzione allora $x$ e $y$ danno luogo ad una rappresentazione primitiva di di $z^l$ come somma di due quadrati.

Per il Corollario 10.28 l'intero Gaussiano $x + iy$ associato ad una scomposizione primitiva in somma di due quadrati di $z^l$ è dato dalla $l-$esima potenza dell'intero Gaussiano $a' + ib'$ associato ad una generica scomposizione primitiva di $z$ in somma di sue quadrati, inoltre risulta essendo $l > 1$ che deve essere $g = 0$ cioè $z$ è dispari. Quindi:

$$x + iy = (a' + ib')^l$$

posto : $a = \max(|a'|, |b'|)$ e $b = \min(|a'|, |b'|)$

abbiamo che:

$x + iy = i^r (a + ib)^l$ con $r = 0,1,2,3$

$MCD(a, b) = 1$ e $a > b > 0$.

Posto

$$\hat{x} + i\hat{y} = (a + ib)^l \qquad (12.3)$$

abbiamo dimostrato il seguente teorema:

*Teorema 12.1*

Per $l > 1$ intero assegnato la generica soluzione primitiva $(x, y, z)$ della (12.1) e definita dalle seguenti formule:

$(x, y, z) \in \{(j_1\hat{x}, j_2\hat{y}, \hat{z}), (j_1\hat{y}, j_2\hat{x}, \hat{z}) \mid j_k = -1,1; k = 1,2\}$ se $l$ dispari.

$(x, y, z) \in \{(j_1\hat{x}, j_2\hat{y}, j_3\hat{z}), (j_1\hat{y}, j_2\hat{x}, j_3\hat{z}) \mid j_k = -1,1; k = 1,2,3\}$ se $l$ pari.

$$\hat{x} = \frac{(a + ib)^l + (a - ib)^l}{2} = \mathrm{Re}\left((a + ib)^l\right) \qquad (12.4)$$

$$\hat{y} = \frac{(a + ib)^l - (a - ib)^l}{2i} = \mathrm{Im}\left((a + ib)^l\right) \qquad (12.5)$$

$\hat{z} = a^2 + b^2$, $MCD(a, b) = 1$, $a$ e $b$ con parità differente, e $a > b > 0$. (12.6)

Dimostrazione

Svolta nei ragionamenti che precedono il teorema. □

*Osservazione 12.2*

Abbiamo di nuovo una generalizzazione delle formule (11.3) che definiscono le terne pitagoriche primitive. Inoltre nel caso di $l = 2n + 1$, cioè numero dispari, esiste una simmetria tra la formula (12.4) la formula (12.5). Infatti osserviamo che $a - ib = i^3(b + ia)$ e quindi che $a + ib = (-i)^3(b - ia)$, sostituendo nella (12.5) si ha che:

$$\hat{y} = \frac{(-1)^{3l} i^{3l}(b - ia) - i^{3l}(b + ia)}{2i} = i^{3(l-1)} \frac{(b + ia)^l + (-1)^{3l-1}(b - ia)^l}{2} = (-i)^{l-1} \frac{(b + ia)^l + (-1)^{l-1}(b - ia)^l}{2}$$

e analoghi passaggi per $\hat{x}$ si possono fare per la (12.4), quindi nel caso di $l$ dispari abbiamo che:

$$\hat{x} = (-1)^{\frac{l+1}{2}} \frac{(b + ia)^l - (b - ia)^l}{2i} = (-1)^{\frac{l+1}{2}} \mathrm{Im}\left((b + ia)^l\right) \qquad (12.7)$$

$$\hat{y} = (-1)^{\frac{l-1}{2}} \frac{(b + ia)^l + (b - ia)^l}{2} = (-1)^{\frac{l-1}{2}} \mathrm{Re}\left((b + ia)^l\right) \qquad (12.8)$$

Quindi le formule e le proprietà che otteniamo per $\hat{x}$ valgono anche per $\hat{y}$ scambiando $a$ con $b$, e viceversa, a parte eventualmente un cambio di segno che dipende dalla classe di resto modulo 4 di $l$.

*Teorema 12.3*

Se l'esponente $l$ è dispari valgono le seguenti proprièta per $\hat{x}$ e $\hat{y}$ definiti nel Teorema 12.1.

1. $\hat{x} = au$ con $MCD(a,u) = 1$
2. $\hat{y} = bv$ con $MCD(b,v) = 1$
3. $u \equiv v \pmod z$ e $MCD(u,v) = 1$
4. $u = V_n(a,b)$, $v = (-1)^n V_n(b,a)$
   con $V_n(X,Y) = \sum_{j=0}^{n} (-1)^j \binom{2n+1}{2j} X^{2(n-j)} Y^{2j}$ e $l = 2n+1$
5. $V_n(X,Y) = \dfrac{(X+Y)^{2n+1} + (X-Y)^{2n+1}}{2X} - 2YR_n(X,Y)$
   dove $R_n(X,Y) = \sum_{j=1,3,\dots}^{\le n} \binom{2n+1}{2j} X^{2(n-j)} Y^{2j}$
6. $(-1)^{n+1} V_n(X,Y) = Z^n + \sum_{j=0}^{n-1} Z^j \left[(Y+iX)^{2(n-j)} + (Y+iX)^{2(n-j)}\right]$
   dove $Z = X^2 + Y^2$.

Dimostrazione

Dalla formula (12.3) usando lo sviluppo del binomio abbiamo:

$$\hat{x} + i\hat{y} = (a+ib)^l = \sum_{j=0}^{l} \binom{l}{k} a^{l-k} (ib)^k$$

separando la parte reale dalla parte immaginaria:

$\hat{x} = \sum_{k=0,2,\dots}^{\le l} \binom{l}{k} a^{l-k} b^k (-1)^{\frac{k}{2}}$ se ora $l = 2n+1$ posto $k = 2j$ abbiamo:

$\hat{x} = a\sum_{j=0}^{n} \binom{2n+1}{2j} a^{2(n-j)} b^{2j} (-1)^j = aV_n(a,b)$ siccome $MCD(a,b) = 1$ e $V_n(a,b) = a\lambda + b^{2n}$ si deve avere $MCD(a, V_n(a,b)) = 1$, tenuto presente quanto detto nella Osservazione 12.2, i punti 1,2 e 4 sono dimostrati.

Ecco di seguito i passaggi che dimostrano il punto 5.:

$$(a+b)^l + (a-b)^l = \sum_{k=0}^{l} \binom{l}{k} a^{l-k} \left(b^k + (-1)^k b^k\right) = 2\sum_{j=0}^{n} \binom{2n+1}{2j} a^{2(n-j)} b^{2j}$$

$$= 2\sum_{j=0,2,\dots}^{\le n} \binom{2n+1}{2j} a^{2(n-j)} b^{2j} + 2\sum_{j=1,3,..}^{\le n} \binom{2n+1}{2j} a^{2(n-j)} b^{2j}$$

$$= 2\sum_{j=1}^{n} \binom{2n+1}{2j} a^{2(n-j)} b^{2j} (-1)^j + 4\sum_{j=1,3,..}^{\le n} \binom{2n+1}{2j} a^{2(n-j)} b^{2j}$$

$$= 2aV_n(a,b) + 4abR_n(a,b)$$

Per la dimostrazione di 6. si tratta di usare la ben nota formula dei prodotti notevoli:

$$\hat{y} = \frac{(a+ib)^l - (a+ib)^l}{2i} = b\sum_{k=0}^{l-1} (a+ib)^k (a-ib)^{l-k-1}$$

segue che:

$$(-1)^n V_n(b,a) = \sum_{k=0}^{l-1} (a+ib)^k (a-ib)^{l-k-1} \qquad (12.9)$$

per $0 < k \le l-k-1$ cioè per $0 < k \le n$ abbiamo che:

$$(a+ib)^k (a-ib)^{l-k-1} = \left(a^2+b^2\right)^k (a-ib)^{2(n-k)} = z^k (a-ib)^{2(n-k)}$$

quindi riscrivendo la (12.9) come segue:

$$\sum_{k=0}^{l-1}(a+ib)^k(a-ib)^{l-k-1}=(a-ib)^{l-1}+\sum_{k=1}^{\frac{l-3}{2}}(a+ib)^k(a-ib)^{l-k-1}+\left(a^2+b^2\right)^{\frac{l-1}{2}}+$$
$$\sum_{k=\frac{l+1}{2}}^{l-2}(a+ib)^j(a-ib)^{l-k-1}+(a+ib)^{l-1}$$

Se nella seconda sommatoria dopo il segno di uguaglianza si pone $k'=l-k-1$ abbiamo che $0<k'\leq n-1$. Segue quindi sostituendo gli indici $k$, $k'$ con l'unico indice $j$ che vale la seguente:

$$(-1)^n V_n(b,a)=z^n+\sum_{j=0}^{n-1}z^j\left[(a+ib)^{2(n-j)}+(a-ib)^{2(n-j)}\right]$$

segue che la 6. risulta dimostrata. Dato che:

$$(a+ib)^{2n}+(a-ib)^{2n}=\left((-i)^3\right)^{2n}(b-ia)^{2n}+\left(i^3\right)^{2n}(b+ia)^{2n}=(-1)^n\left[(b+ia)^{2n}+(b-ia)^{2n}\right]$$

segue:

$$u-v=V_n(a,b)-(-1)^n V_n(b,a)=\lambda z$$

che dimostra la 3.. □

Osservazione 12.4
Dopo il lavoro di Wiles sappiamo che l'ultimo teorema di Fermat è vero per tutti gli esponenti, e possiamo dedurre altre utili corollari dalle formule del Teorema 12.3. Infatti sappiamo in particolare che $\hat{x}$ e $\hat{y}$ non possono essere potenze $l-$esime di numeri interi coprimi. Da questo di deduce facilmente che non esistono due interi coprimi $a_0$ e $b_0$ tali che ${u_0}^l=V_n\left({a_0}^l,{b_0}^l\right)$ e ${v_0}^l=(-1)^n V_n\left({b_0}^l,{a_0}^l\right)$. Ovviamenente considerazioni del genere sono enormemente più facili di una della dimostrazione dell'ultimo teorema di Fermat.

## 13. *L'equazione Diofantea* $X^2 + Y^2 + Z^2 = W^2$ *in* $\mathbb{Z}$

Ci occupiamo di determinare tutte le quaterne di interi $(x, y, z, w)$ soluzioni primitive della equazione diofantea:

$$X^2 + Y^2 + Z^2 = W^2 \quad (13.1)$$

In tal caso soluzione primitiva significa che:

$$MCD(x, y, z, w) = 1 \quad (13.2)$$

Risulta relativamente semplice verificare che se $x^2 + y^2 + z^2 = w^2$:

$$MCD(x, y, z, w) = 1 \Leftrightarrow MCD(y, z, w) = 1 \Leftrightarrow MCD(x, z, w) = 1 \Leftrightarrow MCD(x, y, w) = 1 \Leftrightarrow MCD(x, y, z) = 1 \quad (13.3)$$

*Teorema 13.1-* ***Carmichael [DA1] § 10.***

La generica soluzione $(x, y, z, w)$ della (13.1) tale che $MCD(x, y, z) = 1$ e $w > 0$, eventualmente permutando tra loro e cambiando segno a $x$, $y$, $z$, è definita dalle seguenti formule parametriche:

$$\begin{cases} x = 2(mn - uv) \\ y = m^2 - n^2 - u^2 + v^2 \\ z = 2(mu + nv) \\ w = m^2 + n^2 + u^2 + v^2 \end{cases} \quad (13.4)$$

$$MCD(m, n, u, v) = 1 \quad (13.5)$$

$m, v$ oppure $n, u$ con differente parità.

Dimostrazione

Sia $(x, y, z, w)$ che soddisfa la (13.1) con $MCD(x, y, z) = 1$ e $w > 0$, allora uno ed uno solo tra $x$, $y$ e $z$ può essere dispari, altrimenti risulterebbe $w^2 \equiv 3 \pmod 4$ o $w^2 \equiv 2 \pmod 4$ e entrambi i casi sono impossibili. Permutando se necessario tra loro $x$, $y$, $z$, possiamo assumere:

$$x = 2x' \text{ e } z = 2z'. \quad (13.6)$$

quindi si ha che:

$$4\left(x'^2 + z'^2\right) = w^2 - y^2$$

essendo $w$ e $y$ dispari uno solo tra $w - y$ e $w + y$ è divisibile per $4$ ed entrambi per $2$:

$$\frac{(w + y)}{2}\frac{(w - y)}{2} = x'^2 + z'^2$$

Poniamo:

$d = MCD(x', z')$, $x'' = x'/d$, $z'' = z'/d$.

$$f = \frac{w - y}{2},\ g = \frac{w + y}{2} \quad (13.7)$$

$h = MCD(f, g)$.

Risulta che $MCD(h, d) = 1$, altrimenti si contraddice l'ipotesi $MCD(x, y, z) = 1$. Infatti $f$ e $g$ non possono avere fattori comuni del tipo $3 \pmod 4$ altrimenti questi sarebbero fattori comuni di $w$ e $y$, e dovrebbero dividere $d$ e quindi $x'$ e $z'$. Quindi tutti gli eventuali fattori del tipo $3 \pmod 4$ del numero $x'^2 + z'^2$ si distribuiscono in maniera disgiunta nei due numeri $f$ e $g$ e hanno esponenti pari in essi. Pertanto entrambi i numeri sono scrivibile some somma di due quadrati.

Poniamo $f = (m + iv)(m - iv)$, $g = (n + iu)(n - iu)$ allora abbiamo che:

$$x' + iz' = (m + iv)(n + iu) \quad (13.8)$$

Dal fatto che $f$ e $g$ non possono avere fattori comuni del tipo $3 \pmod 4$ segue che:

$$MCD(MCD(m, v), MCD(n, u)) = MCD(m, n, u, v) = 1$$

Essendo uno ed uno solo tra $f$ e $g$ pari, segue che $m, v$ oppure $n, u$ con differente parità.

Applicando le (13.6),(13.7),(13.8) si ottengono le (13.4). □

*Osservazione 13.2*

Come caso particolare abbiamo di nuovo ritrovato le formule delle terne pigatoriche primitive infatti se $z = 0$, allora $f$ e $g$ sono quadrati e hanno differente parità, ponendo ad esempio $u = v = 0$ ritroviamo le formule (11.3). Un'altra analogia delle quadruple con le triple pitagoriche la si trova se ci poniamo la domanda "quante sono le quadruple primitive per un assegnato $w$?", la cui risposta sta nel numero di rappresentazioni del numero $w$ come somma di quattro quadrati. Nel caso delle triple la risposta stava nel numero di rappresentazioni di $z$ come somma di due quadrati. Lagrange ha dimostrato che ogni numero è scrivibile come somma di quattro quadrati, quindi a differenza delle terne pitagoriche per cui $z$ deve rispettare le condizioni del Teorema 9.10, nel caso delle quadruple pitagoriche $w$ può essere un numero positivo qualsiasi.

## *Bibliografia*

*Testi on-line gratuiti:*

http://us.geocities.com/alex_stef/mylist.html
sotto la voci "Number Theory" e "Algebra, Algebraic Geometry".

http://www.mat.uniroma3.it/ntheory/number_theory.html

http://www.gotmath.com/notes.html

http://historical.library.cornell.edu/math/

*Testi da acquistare o consultare in biblioteca:*